\documentclass[11pt,a4paper]{article}
\usepackage{setspace, hyperref}
\usepackage{amsfonts}
\usepackage{epsfig}
\usepackage{latexsym}
\usepackage{amsthm}
\usepackage{enumitem}
\usepackage{amsmath,stmaryrd}
\usepackage{amssymb}
\newenvironment{keyword}
{\par\noindent\textbf{Keywords:}\ }
{\par}
\newenvironment{msc}
{\par\noindent\textbf{MSC~2020:}\ }
{\par}
\usepackage{array}
\usepackage{tikz}
\usepackage{circuitikz}
\usepackage{tikz-cd}
\usetikzlibrary{shadings,patterns}
\usetikzlibrary{calc,decorations.pathmorphing}
\usepackage{graphicx,graphics}
\newtheorem{theorem}{Theorem}[section]
\newtheorem{lemma}[theorem]{Lemma}
\newtheorem{problem}[theorem]{Problem}
\newtheorem{proposition}[theorem]{Proposition}
\newtheorem{example}[theorem]{Example}
\newtheorem{corollary}[theorem]{Corollary}
\newtheorem{definition}[theorem]{Definition}

\newtheorem{observation}[theorem]{Observation}
\newtheorem{conjecture}[theorem]{Conjecture}
\newtheorem{result}[theorem]{Result}

\def\red#1{{\color{red}#1}}

\def\deg{\mbox{\rm deg}}
\newcommand{\ZC}{Z_C}
\newcommand{\lo}{\operatorname{lo}}
\newcommand{\U}{\mathcal{U}}
\newcommand{\Sk}{\mathcal{S}}
\title{On the connected zero forcing polynomial of finite trees}
\author{Sanjay Mukherjee\thanks{Harish-Chandra Research Institute, A CI of Homi Bhabha National Institute, Chhatnag Road, Jhunsi, Prayagraj - 211019, India. Email: sanjay.math1729@gmail.com}
\and
Papi Ray\thanks{Indian Institute of Science Education and Research Mohali, 140306, India. Email: popiroy93@gmail.com}}

\begin{document}

\maketitle

\begin{abstract}
 Given a graph $\Gamma$ with vertex set $V(\Gamma)$ and a set of colored vertices S, the zero forcing rule is an iterative graph coloring process which states that a colored vertex with exactly one uncolored neighbor forces that neighbor to become colored. If repeated application of this rule eventually colors every vertex of $\Gamma$, the set $S$ is called a \emph{zero forcing set}. Now, if the induced subgraph $\Gamma\vert_S$ is connected, then $S$ is called a \textit{connected zero forcing set}. The minimum cardinality of such a set is the \emph{connected zero forcing number} $Z_C(\Gamma)$. A \emph{connected zero forcing polynomial} of a finite graph $\Gamma$ is a monic polynomial $Z_C^{\Gamma}(x) \;=\; \sum_{k} a_k x^k$, where $a_k$ counts the number of connected zero forcing sets of $\Gamma$ of size $k$, which we denote by $Z_C(\Gamma,k)$. In this paper, we characterize the connected zero forcing polynomial of any finite tree $T$. After that, we prove that the connected zero forcing polynomial $Z_C^T(x)$ of any finite tree $T$ is always log-concave, and thus it is unimodal as well, thereby resolving an open question posed by Newton and Brimkov in [T. Newton, B. Brimkov, Zero forcing is not unimodal.~\emph{arXiv:2611.11352}]. We then establish sharp lower bounds for each coefficient of the connected zero forcing polynomial, over all trees on $n$ vertices with $p$ many pendant vertices. We also identify a family of trees, which we call uniform trees, and we prove that for any tree on $n$ vertices and $p$ many pendant vertices, there exists a uniform tree on $n$ vertices and $p$ many pendant vertices, whose connected zero forcing polynomial bounds as an upper bound, the connected zero forcing polynomial of the tree coefficient wise. Again, for $2\leq p\leq 4$, we establish a  sharp upper bound for each coefficient of the connected zero forcing polynomial, over all trees on $n$ vertices with $p$ many pendant vertices. At the end we consider a few interesting examples which support the fact that for $p\geq 5$, even if $n$ and $p$ are both fixed, in general there is no single tree whose connected zero forcing polynomial dominates coefficient by coefficient the polynomials of every other tree with $n$ vertices and $p$ many pendant vertices.
 
\end{abstract}

        \begin{msc}
05A15, 05A20, 05C05, 05C31.
        \end{msc}
\section{Introduction}

The zero forcing process on graphs, introduced \cite{AIMGroup} originally in connection with the minimum rank problem for symmetric matrices whose off-diagonal pattern is described by a graph, has grown into a vibrant topic of research in its own right. Given a graph $\Gamma$ with vertex set $V(\Gamma)$ and an initial set $S$ of ``colored'' vertices, the color change rule states that a colored vertex with exactly one uncolored neighbor forces that neighbor to become colored. If repeated application of this rule eventually colors every vertex of $\Gamma$, the set $S$ is called a \emph{zero forcing set}, and the minimum cardinality of such a set is the \emph{zero forcing number} $Z(\Gamma)$. Beyond its original role in bounding the maximum nullity of matrices associated with a graph, zero forcing has since found application in controlling networks of coupled quantum spins, in power domination models for monitoring an electrical grid with phasor measurement units, in fault diagnosis of multiprocessor networks, in modeling the diffusion of information and influence in social networks, where forcing a vertex is a proxy for a status propagating once all but one neighbor's status is already determined.

A natural refinement of the process, motivated by settings where the monitoring or controlling agents must themselves form a connected substructure, leads to the notion of connected zero forcing. For a graph $\Gamma$, a subset $S$ of $V(\Gamma)$ is a \emph{connected zero forcing set} if $S$ is a zero forcing set and the induced subgraph $\Gamma\vert_S$ is connected. The minimum cardinality of such a set is the \emph{connected zero forcing number} $Z_C(\Gamma)$. This restriction reflects a genuine constraint in many of the systems zero forcing models: in PMU deployment, for instance, the dominant cost is often the communication infrastructure and maintenance overhead rather than the devices themselves, so a utility benefits from clustering its measurement units into a single connected region, a concern that has motivated a dedicated line of work on connected power domination \cite{brimkovmikesellsmith2019}. A similar principle holds for contagion and influence models such as target set selection, which typically assume the initially active vertices originate from a single, connected source \cite{ackermanbenzwiwolfovitz2010}, making $Z_C(\Gamma)$ the more faithful parameter in that setting.

Beyond the zero forcing number itself, the collection of all zero forcing sets of a graph has been studied from several enumerative and dynamical angles. There is a close relation between Restricted power domination and zero forcing number \cite{bozem}. In \cite{bur}, Burgarth et at investigated how the graph theoretic concept of a zero forcing set impacts the controllability property. In \cite{curt} and \cite{english}, the authors worked on zero forcing on random graphs. In \cite{group} and \cite{huang}, the authors used the zero forcing number to bound the minimum rank or  maximum nullity of a graph. In \cite{goldberg}, Goldberg and Berman provided a bound for the maximum nullity of a matrix with a given sign pattern. The \emph{propagation time} of a zero forcing set records the number of rounds of the color change rule needed to force the entire graph, and one natural question is how large or small this quantity can be over all minimum zero forcing sets \cite{berliner, hogbenkingsleymeyeretal2012, warnberg}. A related parameter, \emph{throttling}, optimizes the trade-off between the size of a zero forcing set and its propagation time by minimizing their sum over all zero forcing sets of the graph \cite{butleryoung2013, carlson}. Both perspectives treat the full family of zero forcing sets, rather than a single minimum-cardinality representative, as the object of interest. Based on this observation, in \cite{boyer} Boyer et al.~introduced the notion of the \emph{zero forcing polynomial} of a finite graph as a monic polynomial, where the coefficient of the $k$-th exponent in the polynomial enumerates the number of zero forcing sets with $k$ vertices in that graph. Clearly, the exponent corresponding to the smallest non-zero coefficient in that polynomial denotes the zero forcing number of the graph.

An analogous refinement exists for the connected variant. In \cite{bd}, Brimkov and Davila discussed various applications of the connected zero forcing number. Furthermore, in \cite{bh}, Brimkov and Hicks proved that determining whether a given initial colored set can successfully complete the zero forcing process remains NP-hard, even when the initial colored set induces a connected subgraph.  In \cite{khosravirashidisheikhhosseni2020} Khosravi, Rashidi, and Sheikhhosseni introduced the \emph{connected propagation time} of a graph and characterized the graphs attaining its extreme values, showing that questions of propagation speed remain meaningful once the initial set is required to stay connected throughout the process. A connected counterpart to throttling, by contrast, appears comparatively unexplored. Taken together, these results indicate that a fuller accounting of connected zero forcing sets beyond just their minimum cardinality, is a natural and largely open direction, one that the connected zero forcing polynomial is designed to address. A \emph{connected zero forcing polynomial} of a finite graph $\Gamma$ is a monic polynomial $Z_C^{\Gamma}(x) \;=\; \sum_{k} a_k x^k$, where $a_k$ counts the number of connected zero forcing sets of $\Gamma$ of size $k$. While the ordinary zero forcing polynomial has by now been studied for several graph families, the connected variant carries genuinely different combinatorics, since connectivity of $S$ interacts heavily with the structure of the underlying graph, and the enumeration problem is correspondingly different even for simple graph class, such as trees. 

In this article, we undertake a systematic study of the connected zero forcing polynomial of trees. We begin by giving an enumeration of $Z_C^T(x)$ for a tree $T$. Building on this, we turn to the question of how the coefficients of $Z_C^T(x)$ behave when the number of vertices $n$ and the number of pendant vertices $p$ of $T$ are held fixed. Having settled the structure of the connected zero forcing  polynomial, we turn to a finer structural property of these polynomials, namely \emph{log-concavity}. After that, we establish sharp lower bounds for each coefficient of the connected zero forcing polynomial, over all trees on $n$ vertices with $p$ many pendant vertices, together with a characterization of the extremal trees attaining these bounds. We then take up the more delicate question of an upper bound. Here the situation is fundamentally different: we show that if $n$ and $p$ are both fixed, then for $p<5$ there is a unique maximal tree but for $p\geq 5$ in general no single tree whose connected zero forcing polynomial dominates coefficient by coefficient the polynomials of every other tree with $n$ vertices and $p$ many pendant vertices. Despite this obstruction, we identify a family of trees, which we call uniform trees, characterized by having their pendant vertices distributed as evenly as possible among the branches at each vertex of degree at least three. We prove that for any tree on $n$ vertices and $p$ many pendant vertices, there exists a uniform tree on $n$ vertices and $p$ many pendant vertices, whose connected zero forcing polynomial bounds the connected zero forcing polynomial of the tree coefficient wise, from above. This gives, in place of a single extremal tree, an explicit and structurally natural family that plays the extremal role for the entire class.

\section{Preliminary}

Before going into the main results, let us first define some preliminaries. A polynomial $f(x)=\sum_j f_jx^j$ with nonnegative coefficients is called \emph{log-concave} if $f_j^2\ge f_{j-1}f_{j+1}$ for all $j$. $f(x)$ is said to have \emph{no internal zeros} if $f_af_b>0$ with $a<b$ implies $f_j>0$ for all $a\leq j\leq b$. For an integer $h\geq 1$, we denote $[h]=1+x+\dots+x^{h-1}$, and $[0]=0$. For $a\in\mathbb Z^k$, define $|a|:=a_1+\dots+a_k$. A nonzero polynomial $A(x)=\sum_{j=0}^{d}A_jx^j$ of degree $d$ is \emph{symmetric} if $A_j=A_{d-j}$ for $0\leq j\leq d$, and \emph{unimodal} if $A_0\leq A_1\le\dots\leq A_p\ge A_{p+1}\ge\dots\ge A_d$ for some $p$. (Andrews \cite{An} calls symmetric polynomials \emph{reciprocal}.) A symmetric polynomial has $A_0=A_d\neq0$. For instance, $[L]$ with $L\ge 1$ is symmetric and unimodal of degree $L-1$, whereas $x^a[L]$ with $a\ge 1$ is \emph{not} symmetric. Accordingly, we never call a polynomial with zero constant term symmetric. Such a polynomial is always written as an explicit power of $x$ times a polynomial with nonzero constant term. Now, the following well known result is useful. Also, the following fact is easy to observe.

\begin{lemma}[\cite{An}; \cite{St}, Prop.~1] \label{Lemma: Product symmetric unimodal}
If $A$ and $B$ are symmetric unimodal polynomials with nonnegative coefficients, then $AB$ is symmetric and unimodal.
\end{lemma}

Let $\Gamma$ be a graph with vertex set $V(\Gamma)$ and edge set $E(\Gamma)$. A vertex $v$ of $\Gamma$ is called a \textit{pendant vertex} if $d(v)=1$. A path of $\Gamma$ is said to be a \textit{pendant path}, if at least one of the extreme vertices of the path is a pendant vertex. A vertex $v$ is called a \textit{branch vertex}, if $d(v)\geq 3$. For a subgraph $S$ of $\Gamma$ and a vertex $u\in V(\Gamma)$, the \textit{external neighborhood} of $u$ in $\Gamma$, denoted $N^{\Gamma}_S(u)$, is defined to be the collection of all vertices in $V(\Gamma)\setminus V(S)$ which are adjacent to $u$ in $\Gamma$. We define $N^{\Gamma}_S[u]:=N^{\Gamma}_S(u)\cup\{u\}$. The \textit{external degree} of $u$, denoted by $d_S(u)$ is defined as $d_S(u)=\lvert N^{\Gamma}_S(u)\rvert$. 

A tree is called a \emph{starlike graph}, if it has exactly one vertex with degree more than two. It is easy to see that a tree has exactly one leaf if and only if it is a trivial graph, exactly two pendant vertices if and only if it is a path. Further, any tree with only three pendant vertices is a starlike graph. For a starlike graph $\Gamma$ with central vertex $u$ and pendant vertices $v_1,v_2,\ldots,v_k$, let $h_i={\rm dist}(u,v_i)$ for $i\in [k]$. Then the starlike graph $\Gamma$ is denoted by $\mathcal{S}_{h_1,h_2,\ldots,h_k}$. 

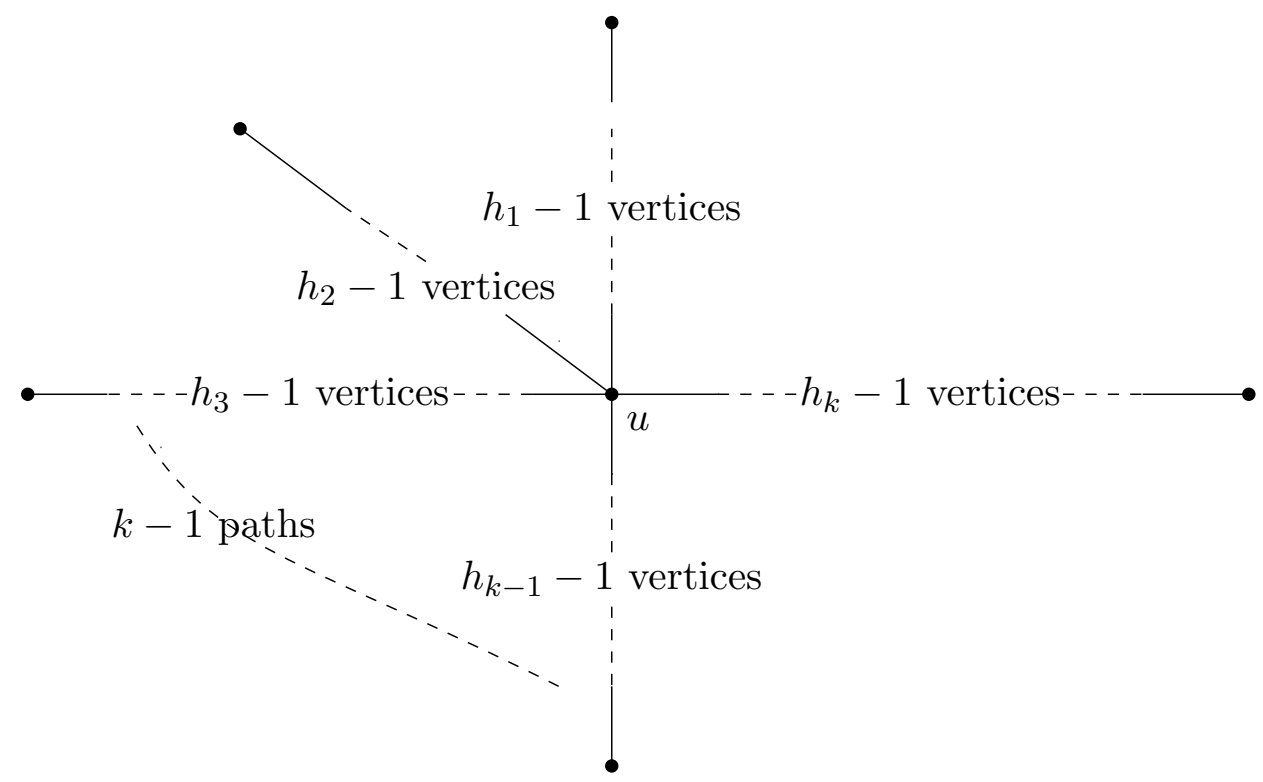
\begin{figure}[!ht]
\centering
\resizebox{0.5\textwidth}{!}{%
\begin{circuitikz}
\tikzstyle{every node}=[font=\normalsize]
\node at (5.5,7.5) [circ] {};
\draw (5.5,7.5) to[short] (4.75,7.5);
\draw (5.5,7.5) to[short] (5.5,6.75);
\draw (5.5,7.5) to[short] (5.5,8.25);
\draw [dashed] (5.5,8.25) -- (5.5,9);
\draw [dashed] (4.75,7.5) -- (4,7.5);
\draw [dashed] (5.5,6.75) -- (5.5,6.75);
\draw [dashed] (5.5,6.75) -- (5.5,6);
\draw [dashed] (1.25,7) -- (1.25,7);
\draw [short] (5.5,7.5) -- (6.5,7.5);
\draw [dashed] (6.5,7.5) -- (6.5,7.5);
\draw [dashed] (6.5,7.5) -- (7.25,7.5);
\node [font=\normalsize] at (8.5,7.5) {$h_k-1$ vertices};
\node [font=\normalsize] at (5.5,9.25) {$h_1-1$ vertices};
\node [font=\normalsize] at (3.75,8.5) {$h_2-1$ vertices};
\draw [short] (5.5,7.5) -- (5.5,7.5);
\draw [short] (5,8) -- (5,8);
\draw [short] (5.5,7.5) -- (4.5,8.25);
\draw [dashed] (3.75,8.75) -- (3,9.25);
\draw [dashed] (5.5,9.5) -- (5.5,10);
\node [font=\normalsize] at (2.75,7.5) {$h_3-1$ vertices};
\draw [dashed] (1.5,7.5) -- (0.75,7.5);
\draw [dashed] (5.5,5.5) -- (5.5,4.75);
\draw [dashed] (9.75,7.5) -- (10.5,7.5);
\draw [short] (10.5,7.5) -- (11.5,7.5);
\draw [short] (5.5,10.25) -- (5.5,11);
\draw [short] (3,9.25) -- (2,10);
\draw [short] (0.75,7.5) -- (0,7.5);
\draw [short] (5.5,4.75) -- (5.5,4);
\node at (11.5,7.5) [circ] {};
\node at (5.5,11) [circ] {};
\node at (2,10) [circ] {};
\node at (0,7.5) [circ] {};
\node at (5.5,4) [circ] {};
\node [font=\normalsize] at (5.5,5.75) {$h_{k-1}-1$ vertices};
\draw [dashed] (5,4.75) .. controls (2.5,6) and (1.75,6) .. (1,7.25);
\node [font=\normalsize] at (1.75,6.25) {$k-1$ paths};
\node [font=\normalsize] at (5.75,7.25) {$u$};
\end{circuitikz}
}%
\caption{The starlike graph $\mathcal{S}_{h_1,h_2,\ldots,h_k}$.}
\label{fig:my_label}
\end{figure}

Let $p\geq 4$ and $n\geq p+2$. Choose a partition $p=k+(p-k)$, where $2\leq k\leq p-2$. Then the graph $\mathcal{H}^n_{k,p-k}$ is the tree with $n$ vertices and $p$ many pendant vertices, such that $\mathcal{H}^n_{k,p-k}$ has exactly two vertices $u_1,u_2$ with degree $\geq 3$, $d(u_1)=k+1$, $d(u_2)=p-k+1$, and every pendant vertex is adjacent to either $u_1$ or $u_2$. 

\begin{figure}[!ht]
\centering
\resizebox{0.5\textwidth}{!}{%
\begin{circuitikz}
\tikzstyle{every node}=[font=\normalsize]
\node at (5,7.75) [circ] {};
\draw (5,7.75) to[short] (5,6.5);
\draw (5,7.75) to[short] (6,7.75);
\draw [dashed] (6.25,7.75) -- (7.75,7.75);
\draw [short] (7.75,7.75) -- (9,7.75);
\draw [short] (6,7.75) -- (6,7.75);
\draw [short] (6,7.75) -- (6,7.75);
\draw [short] (6,7.75) -- (6,7.75);
\node at (9,7.75) [circ] {};
\node at (5,6.5) [circ] {};
\draw (9,7.75) to[short] (9,9);
\draw (9,7.75) to[short] (10.25,7.75);
\draw (9,7.75) to[short] (9,6.5);
\node at (9,6.5) [circ] {};
\draw (5,7.75) to[short] (5,9);
\node at (5,9) [circ] {};
\node at (9,9) [circ] {};
\node at (10.25,7.75) [circ] {};
\draw (5,7.75) to[short] (4,8.75);
\draw (9,7.75) to[short] (10,8.75);
\node at (4,8.75) [circ] {};
\node at (10,8.75) [circ] {};
\draw [dashed] (4.75,6.75) .. controls (3.75,7.25) and (3.75,7.75) .. (4.25,8.5);
\draw [dashed] (9.25,6.75) .. controls (10,7) and (10.25,7.25) .. (10,7.75);
\node [font=\normalsize] at (10,7.25) {$k$ pendant vertices};
\node [font=\normalsize] at (3.75,7.25) {$p-k$ many pendant vertices};
\node [font=\normalsize] at (8.75,7.5) {$u_1$};
\node [font=\normalsize] at (5.25,7.5) {$u_2$};
\end{circuitikz}
}%
\caption{The graph $\mathcal{H}^n_{k,p-k}$.}
\label{fig:the h graph}
\end{figure}
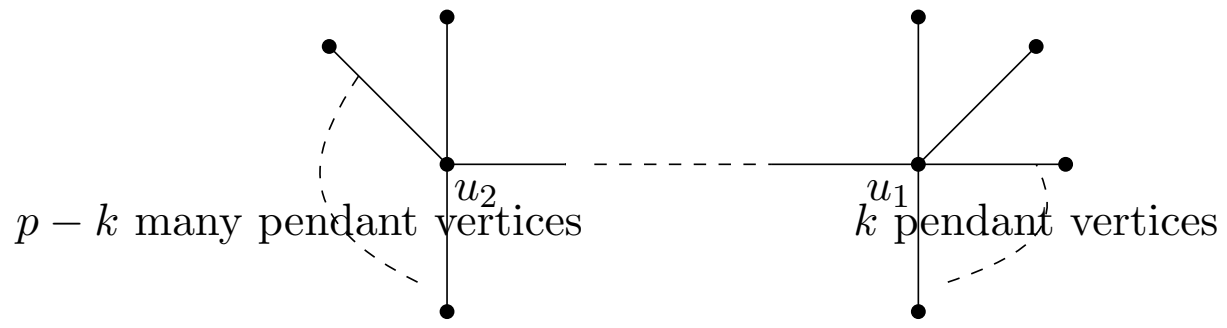
\begin{definition}
For a finite tree $T$, we define the \textit{stem} of $T$ to be the minimal subtree of $T$ containing all the branch vertices of $T$. We denote the stem of $T$ by ${\rm st}(T)$.     
\end{definition}
Note that the stem of a finite tree is the union of all paths between the branch vertices of that tree. Hence, for a finite tree $T$, ${\rm st}(T)$ is unique. 

For a connected graph $\Gamma$ and a subset $S\subseteq V(\Gamma)$, let $\mathfrak{Z}_C^S(\Gamma;i)$ denote the collection of all connected zero forcing sets of $\Gamma$ of size $i$ containing $S$. When $S=\emptyset$, i.e.~the collection of all connected zero forcing sets of $\Gamma$ of size $i$ is denoted by $\mathfrak{Z}_C(\Gamma;i)$. Further, we denote $\lvert\mathfrak{Z}_C^S(\Gamma;i)\rvert={Z}_C^S(\Gamma;i)$ and $\lvert\mathfrak{Z}_C(\Gamma;i)\rvert={Z}_C(\Gamma;i)$. We define $Z_C^{\Gamma,S}(x):= \sum_{i=0}^{\lvert V(\Gamma) \rvert}{Z}_C^S(\Gamma;i)x^i$ to be the polynomial having coefficients as the number of connected zero forcing sets of the size of the corresponding exponent, containing the set $S$. Also, when $S$ contain a single vertex $v$, we write $S$ as $v$ in place of $\{v\}$.

Let $P(\Gamma)$ denote the set of pendant vertices of $\Gamma$ and $p(\Gamma):=\lvert P(\Gamma)\rvert$. We denote by $\mathfrak{T}_{n,p}$ the collection of all trees with $n$ vertices and $p$ many pendant vertices. A pendant vertex $l$ of a tree $T$ is said to be a \textit{seed} of a vertex $v\in V({\rm st}(T))$, if the pendant path between $v$ and $l$ contains only the vertex $v$ from $v\in V({\rm st}(T))$. We denote by the set $L_v$, the set of all seeds of $v$ in $T$. If $\lvert L_v\rvert=1$ for some vertex $v\in {\rm st}(T)$, we denote that only vertex by $l_v$. Note that a leaf in $T$ has a unique seeding vertex in ${\rm st}(T)$. By $\mathfrak{d}(l)$, we denote the distance of the leaf $l$ from its seeding vertex. For a tree $T$, we define the following.
\begin{center}
    $S_1^T:=\{v\in V({\rm st}(T))~\vert~d_{{\rm st}(T)}(v)=1\}$.\\
    $S_2^T:=\{v\in V({\rm st}(T))~\vert~d_{{\rm st}(T)}(v)=2\}$.\\
    $S_3^T:=\{v\in V({\rm st}(T))~\vert~d_{{\rm st}(T)}(v)\geq 3\}$.
\end{center}
Define $S^T:=S_1^T\cup S_2^T \cup S_3^T$. For a vertex $u\in S^T$, let $T_u$ denote the minimal subtree of $T$ containing $u$ and all the seeds of $u$. It is worth noticing that for any vertex $v\in S_1^T\cup S_2^T$, $T_v$ is a path, and for any vertex $w\in S_3^T$, $T_w$ is a starlike graph.

\begin{observation}\label{distance tree size observation}
$\lvert V({\rm st}(T))\rvert+\underset{l\in P(T)}{\sum}\mathfrak{d}(l)=n$.    
\end{observation}

\begin{observation}\label{leaf seed lemma}
Every leaf of $T$ has a unique seeding vertex. Also, for every vertex $v$ of ${\rm st}(T)$, we have $\lvert L_v\rvert=d_{{\rm st}(T)}(v)$.
\end{observation}

As a corollary of the above, we get the following.
\begin{corollary}\label{leaf external degree corollary}
Let $T\in \mathfrak{T}_{n,p}$. Then $\underset{v\in V({\rm st}(T))}{\sum}d_{{\rm st}(T)}(v)=p$.    
\end{corollary}

\section{The connected zero forcing polynomial of trees}

Trees give the connected variant a particularly clean structure, since any connected zero forcing set of a tree $T$ induces a subtree of $T$, and forcing all of $T$ reduces to identifying a small connected ``root'' from which the coloring spreads along every branch. This makes trees a natural setting for moving past the connected zero forcing number itself and asking a more refined question---not just how small a connected zero forcing set can be, but how many such sets of each size a tree admits---which is exactly what the connected zero forcing polynomial is built to capture. 

The following result is a consequence of \cite[Observation 2]{li}.
\begin{lemma}\label{zero forcing set containing stem}
Every connected zero forcing set of a tree $T$ contains $V({\rm st}(T))$.    
\end{lemma}
 
Throughout the rest of this section, we compute the connected zero forcing polynomial of trees. Before getting into the main results, let us first observe the following facts.
\begin{lemma}\label{connected_zero_forcing_path}
Let $P$ be a path of length $n$, with two endpoints $v_1$ and $v_2$, and an intermediate point $w$ at a distance of $l$ from $v_1$, with $0<l<n$. Then $Z_C^{P,v_1}(x)=Z_C^{P,v_2}(x)=x(1+x\frac{x^n-1}{x-1})$, and $Z_C^{P,w}(x)=x^2\frac{(x^k-1)(x^{n-k}-1)}{(x-1)^2}(x+(x-1)(\frac{1}{x^l-1}+\frac{1}{x^{n-l}-1}))$.   
\end{lemma}
\begin{proof}
Since, $v_1$ (equivalently $v_2$) is a terminal point of $P$, for any $k\in\{1,2,\ldots,n\}$ the only connected zero forcing set of $P$ containing $v_1$ (equivalently $v_2$) is the subpath of $P$ starting at $v_1$ (equivalently $v_2$), with length $k$. So, we have $Z_C^{P,v_1}(x)=Z_C^{P,v_2}(x)= x+x^2+\cdots+x^{n+1}=x(1+x\frac{x^n-1}{x-1})$. 

Now, consider the set of all connected zero forcing sets of $P$, containing $w$. It is easy to see that a connected zero forcing set of $P$, containing $w$ can miss at most one neighbor of $w$. Now, the generating function for the number of connected zero forcing sets of $P$ containing both the neighbors of $w$ is $x(x+x^2+\cdots+x^l)(x+x^2+\cdots+x^{n-l})$, and generating function for the number of connected zero forcing sets of $P$ containing exactly one of the neighbors of $w$ is $x(x+x^2+\cdots+x^l)+x(x+x^2+\cdots+x^{n-l})$. Thus, we have $Z_C^{P,w}(x)=x(x+x^2+\cdots+x^l)(x+x^2+\cdots+x^{n-l})+x(x+x^2+\cdots+x^l)+x(x+x^2+\cdots+x^{n-l})=x^2\frac{(x^k-1)(x^{n-k}-1)}{(x-1)^2}(x+(x-1)(\frac{1}{x^l-1}+\frac{1}{x^{n-l}-1}))$. 
\end{proof}

\begin{lemma}\label{connected_zero forcing_starlike}
Let $T$ be the starlike graph $S_{h_1,h_2,\dots, h_k}$, for $k\geq 3$. Then $Z_C^T(x)=x^k(\underset{i=1}{\overset{k}{\prod}}\frac{x^{h_i}-1}{x-1})(x+\underset{i=1}{\overset{k}{\sum}}\frac{x-1}{x^{h_i}-1})$.    
\end{lemma}
\begin{proof}
Let $v$ be the central vertex of the starlike graph $S_{h_1,h_2,\dots, h_k}$. Again, it can be seen that a connected zero forcing set of $S_{h_1,h_2,\dots, h_k}$ contains $v$, and can miss at most one neighbor of $v$. Now, the generating function for the number of connected zero forcing sets of $S_{h_1,h_2,\dots, h_k}$ containing all the neighbors of $w$ is $x(\prod_{i=1}^k(x+x^2+\cdots + x^i))$, and the generating function for the number of connected zero forcing sets of $S_{h_1,h_2,\dots, h_k}$ containing all but one neighbors of $w$ is $x(\sum_{i=1}^k\prod_{\substack{j=1 \\ j\neq i}}^k(x+x^2+\cdots + x^i))$. Thus, the connected zero forcing polynomial is $x(\prod_{i=1}^k(x+x^2+\cdots + x^i))+x(\sum_{i=1}^k\prod_{\substack{j=1 \\ j\neq i}}^k(x+x^2+\cdots + x^i))$, and a few steps of simple calculations yield that $Z_C^{S_{h_1,h_2,\dots, h_k}}(x)=x^k(\underset{i=1}{\overset{k}{\prod}}\frac{x^{h_i}-1}{x-1})(x+\underset{i=1}{\overset{k}{\sum}}\frac{x-1}{x^{h_i}-1})$.
\end{proof}
Now, we have the following. 

\begin{proposition}\label{Prop:zero_forcing_polynomial_of_trees_structure}
Let $T$ be a finite tree. Then $$Z_C^T(x)=x^{\rvert V({\rm st}(T))\lvert-\lvert S^T\rvert}(\underset{v\in S_1^T\cup S_2^T}{\prod}Z_C^{T_v,v}(x))(\underset{w\in S_3^T}{\prod}Z_C^{T_w}(x)).$$    
\end{proposition}
\begin{proof}
Let $T$ be a finite tree. Then Lemma \ref{zero forcing set containing stem} implies that any connected zero forcing set of $T$ contains $V({\rm st}(T))$. Let $v\in S^T$. Then Figure \ref{fig:The_subtrees} shows that if $v\in S_1^T\cup S_2^T$, then the restriction of any connected zero forcing set of $T$ over $T_v$ is a connected zero forcing set of $T_v$, containing $v$, and if $v\in S_3^T$, then the restriction of any connected zero forcing set of $T$ over $T_v$ is a connected zero forcing set of $T_v$. In addition, it can be seen that taking a collection of a zero forcing set of $T_v$ containing $v$, for all $v\in S_1^T\cup S_2^T$, and a collection of a zero forcing set of $T_w$ for all $w\in S_3^T$, along with all the remaining vertices of ${\rm st}(T)$, forms a connected zero forcing set of $T$. Thus, the connected zero forcing polynomial of $T$ is a product of $Z_C^{T_v,v}(x)$ for all $v\in S_1^T\cup S_2^T$, with $Z_C^{T_w}(x)$ for all $w\in S_3^T$, multiplied by $x^{\rvert V({\rm st}(T))\lvert-\lvert S^T\rvert}$. This concludes the proof.

\begin{figure}
    \centering
\begin{tikzpicture}[scale=0.8,
    x=1cm,
    y=1cm,
    line cap=round,
    line join=round
]

\tikzset{
    graph edge/.style={
        line width=1pt
    },
    graph dashed/.style={
        line width=1pt,
        dash pattern=on 4.5pt off 3.2pt
    },
    boundary/.style={
        black!60,
        line width=0.75pt,
        dash pattern=on 0.8pt off 2.0pt
    },
    filled vertex/.style={
        fill=black,
        draw=black,
        line width=0.8pt
    },
    open vertex/.style={
        fill=white,
        draw=black,
        line width=1.15pt
    }
}


\def\axis{3.23}
\newcommand{\mir}[1]{2*\axis-#1}

\draw[boundary]
  (\axis,2.13)
  .. controls (3.72,2.10) and (4.22,2.28)  .. (4.36,2.68)
  .. controls (4.48,3.28) and (4.58,4.10)  .. (4.42,4.72)
  .. controls (4.28,5.28) and (4.02,5.70)  .. (3.62,5.72)
  .. controls (3.42,5.96) and (3.05,5.96)  .. (\axis,5.94)
  .. controls ({\mir{3.05}},5.96)  and ({\mir{3.42}},5.96) .. ({\mir{3.62}},5.72)
  .. controls ({\mir{4.02}},5.70)  and ({\mir{4.28}},5.28) .. ({\mir{4.42}},4.72)
  .. controls ({\mir{4.58}},4.10)  and ({\mir{4.48}},3.28) .. ({\mir{4.36}},2.68)
  .. controls ({\mir{4.22}},2.28)  and ({\mir{3.72}},2.10) .. (\axis,2.13)
  -- cycle;

\draw[graph edge]
(3.18,2.25) -- (3.18,3.22);

\draw[graph edge]
(3.18,3.22) -- (3.18,3.98);

\draw[graph dashed]
(3.18,3.98) -- (3.18,5.48);

\draw[graph edge]
(3.18,2.25) -- (2.16,1.73);

\draw[graph dashed]
(2.16,1.73) -- (1.22,1.16);

\draw[graph edge]
(3.18,2.25) -- (4.36,1.72);

\draw[graph dashed]
(4.36,1.72) -- (5.48,1.17);


\fill[filled vertex] (3.18,2.25) circle (0.105cm);

\fill[filled vertex] (3.18,3.22) circle (0.105cm);

\draw[open vertex] (3.18,3.98) circle (0.115cm);
\draw[open vertex] (3.18,5.48) circle (0.115cm);

\fill[filled vertex] (2.16,1.73) circle (0.105cm);
\fill[filled vertex] (4.36,1.72) circle (0.105cm);

\node at (3.18,1.85) {$v$};
\node at (3.18,0.87)
    {$v\in S_{1}^{T}$};
\node at (4.46,5.48) {$T_v$};

\def\axis{13.15}
\begin{scope}[shift={(-.25,0)}]
\draw[boundary]
  (13.42,2.22)
  .. controls (14.00,2.12) and (14.52,2.55) .. (14.72,3.25)
  .. controls (14.90,4.02) and (14.90,4.82) .. (14.62,5.45)
  .. controls (14.28,6.00) and (13.60,6.25) .. (\axis,6.25)
  .. controls ({\mir{13.60}},6.25) and ({\mir{14.28}},6.00) .. ({\mir{14.62}},5.45)
  .. controls ({\mir{14.90}},4.82) and ({\mir{14.90}},4.02) .. ({\mir{14.72}},3.25)
  .. controls ({\mir{14.52}},2.55) and ({\mir{14.00}},2.12) .. ({\mir{13.42}},2.22)
  -- cycle;
\end{scope}

\draw[graph edge]
(12.88,2.32) -- (11.90,1.78);

\draw[graph dashed]
(11.90,1.78) -- (10.98,1.28);

\draw[graph edge]
(12.88,2.32) -- (14.10,1.77);

\draw[graph dashed]
(14.10,1.77) -- (15.10,1.37);


\draw[graph edge]
(12.88,2.32) -- (12.68,3.18);

\draw[graph edge]
(12.68,3.18) -- (12.52,3.95);

\draw[graph dashed]
(12.52,3.95) -- (12.15,5.30);


\draw[graph edge]
(12.88,2.32) -- (13.16,3.18);

\draw[graph edge]
(13.16,3.18) -- (13.36,3.96);

\draw[graph dashed]
(13.36,3.96) -- (13.70,5.30);


\fill[filled vertex] (12.88,2.32) circle (0.105cm);

\fill[filled vertex] (11.90,1.78) circle (0.105cm);
\fill[filled vertex] (14.10,1.77) circle (0.105cm);

\draw[open vertex] (12.68,3.18) circle (0.115cm);
\draw[open vertex] (12.52,3.95) circle (0.115cm);
\draw[open vertex] (12.15,5.30) circle (0.115cm);

\fill[filled vertex] (13.16,3.18) circle (0.105cm);
\draw[open vertex] (13.36,3.96) circle (0.115cm);
\draw[open vertex] (13.70,5.30) circle (0.115cm);

\node at (12.88,1.9) {$v$};
\node at (12.88,0.87)
    {$v\in S_{2}^{T}$};
\node at (14.71,5.48) {$T_v$};
\end{tikzpicture}
\end{figure}
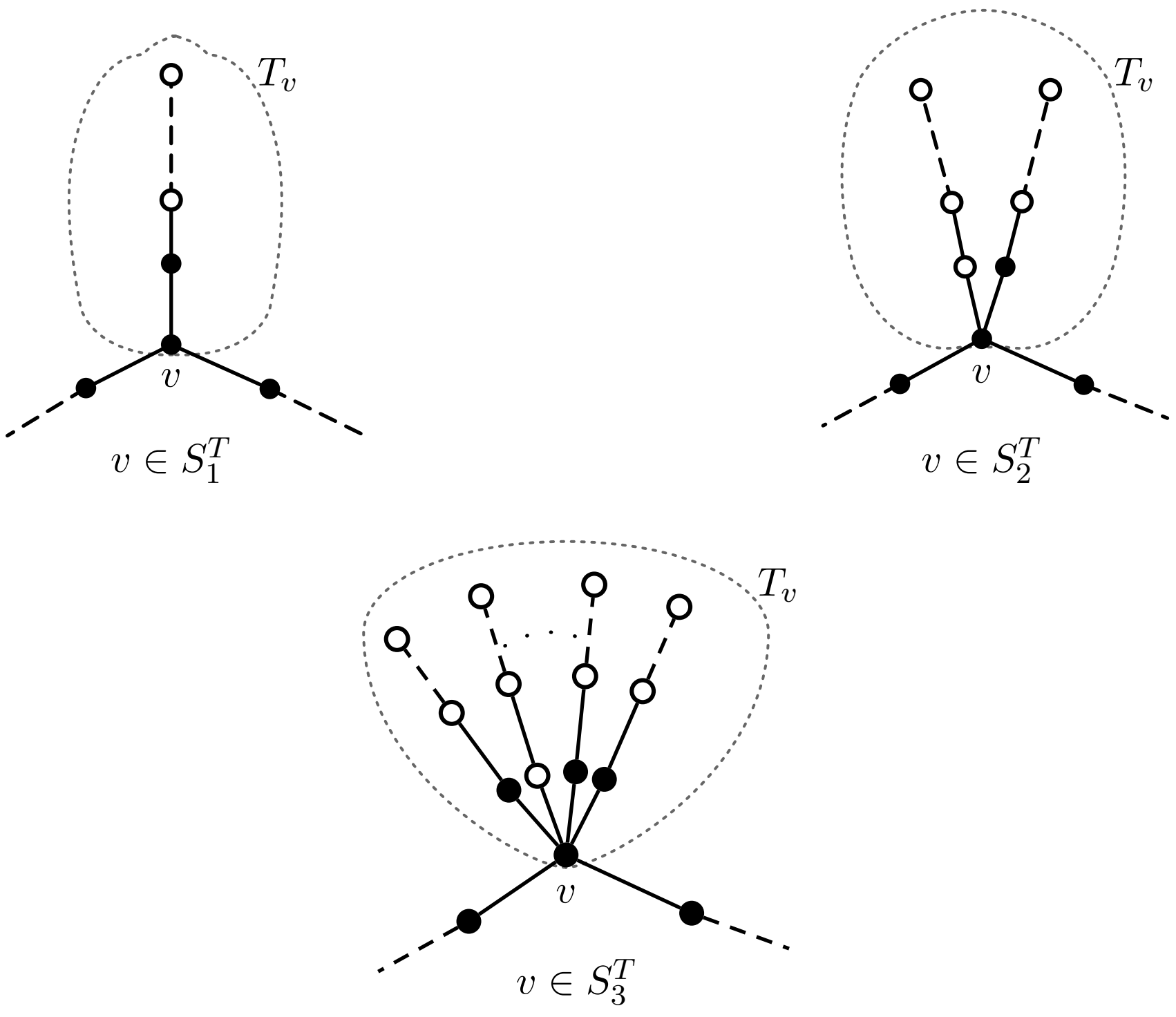
\begin{figure}
    \centering
\begin{tikzpicture}[scale=0.8,
    line width=1pt,
    vtx/.style  = {circle, draw, inner sep=0pt, minimum size=2.1mm},
    blk/.style  = {vtx, line width=0.8pt, fill=black},   
    wht/.style  = {vtx, line width=1.15pt, fill=white},   
    dsh/.style  = {line width=1pt, dash pattern=on 4.5pt off 3.2pt},
    dts/.style  = {line width=1pt, line cap=round, dash pattern=on 0pt off 6.5pt},
    bnd/.style  = {black!60, line width=0.75pt, line cap=round, dash pattern=on 0.8pt off 2.0pt},
  ]

  \def\aA{126.5}\def\aB{107.4}\def\aC{84.3}\def\aD{66.5}

  \pgfmathsetmacro{\bw}{3.4*max(abs(cos(\aA)),abs(cos(\aD)))+0.4}

  \draw[bnd] (0,0)
      .. controls (-0.45,0)        and (-\bw,1.1)      .. (-\bw,2.8)
      .. controls (-\bw,3.4)       and (-0.58*\bw,3.9) .. (0,3.9)
      .. controls ( 0.58*\bw,3.9)  and ( \bw,3.4)      .. ( \bw,2.8)
      .. controls ( \bw,1.1)       and ( 0.45,0)       .. cycle;

  \node[blk] (v) at (0,0.15) {};
  \node[below=0.75pt] at (v.south) {$v$};

  \foreach \i/\ang/\st in {A/\aA/blk, B/\aB/wht, C/\aC/blk, D/\aD/blk}{
    \node[\st] (\i1) at (\ang:1.15) {};
    \node[wht] (\i2) at (\ang:2.3)  {};
    \node[wht] (\i3) at (\ang:3.4)  {};
    \draw (v) -- (\i1) -- (\i2);
    \draw[dsh] (\i2) -- (\i3);
  }

  \draw[dts] (\aB-2:2.75) to[bend left=25] (\aC+2:2.75);

  \node[blk] (E) at (209:1.33) {};  \draw (v) -- (E);  \draw[dsh] (E) -- (209:2.57);
  \node[blk] (F) at (-20:1.6)  {};  \draw (v) -- (F);  \draw[dsh] (F) -- (-20:2.84);

  \node[] at ([yshift=-0.1cm]current bounding box.south) {$v\in S_3^{T}$};
  \node at (2.54,3.38) {$T_v$};
\end{tikzpicture}
 \caption{The subtrees $T_v$.}
    \label{fig:The_subtrees}
\end{figure}
\end{proof}

Using the above results, we now compute the expression for the connected zero forcing polynomial of a tree.

\begin{theorem}\label{Thm: zero forcing polynomial of trees}
Let $T$ be a tree with $n$ vertices and $p$ many pendant vertices. Then $$Z_C^{T}(x)=x^{\rvert V({\rm st}(T))\lvert+p-\lvert S^T\rvert}\frac{\prod_{l\in P(T)}(x^{\mathfrak{d}(l)}-1)}{(x-1)^p}\times\underset{v\in S^T}{\prod}(x+\sum_{l\in L_v}\frac{x-1}{x^{\mathfrak{d}(l)}-1}).$$  
\end{theorem}

\begin{proof}
Let $T$ be a tree with $n$ vertices and $p$ many pendant vertices, with $Z_C(T)=k$. Then Proposition \ref{Prop:zero_forcing_polynomial_of_trees_structure} says that $Z_C^T(x)=x^{\rvert V({\rm st}(T))\lvert-\lvert S^T\rvert}(\underset{u\in S_1^T\cup S_2^T}{\prod}Z_C^{T_v,v}(x))(\underset{w\in S_3^T}{\prod}Z_C^{T_w}(x)).$ So, Lemma \ref{connected_zero_forcing_path} and \ref{connected_zero forcing_starlike} implies
\begin{align*}
Z_C^{T}(x) &=x^{\rvert V({\rm st}(T))\lvert-\lvert S^T\rvert}\underset{u\in S_1^T}{\prod}x(1+x\frac{x^{\mathfrak{d}(l)}-1}{x-1})\times\underset{v\in S_2^T\cup S_3^T}{\prod}(x^{d_{\rm st(T)}(v)}\underset{l\in L_{v}}{\prod}((\frac{x^{\mathfrak{d}(l)}-1}{x-1})(x+\underset{l\in L_{v}}{\sum}(\frac{x-1}{x^{\mathfrak{d}(l)}-1}))))\\
&=x^{\rvert V({\rm st}(T))\lvert+p-\lvert S^T\rvert}\underset{l\in P(T)}{\prod}(\frac{x^{\mathfrak{d}(l)}-1}{x-1})\times\underset{u\in S_1^T}{\prod}(x+\frac{x-1}{x^{\mathfrak{d}(l)}-1})\times\underset{v\in S_2^T\cup S_3^T}{\prod}(x+\underset{l\in L_{v}}{\sum}(\frac{x-1}{x^{\mathfrak{d}(l)}-1}))\\
&=x^{\rvert V({\rm st}(T))\lvert+p-\lvert S^T\rvert}\frac{\prod_{l\in P(T)}(x^{\mathfrak{d}(l)}-1)}{(x-1)^p}\times\underset{v\in S^T}{\prod}(x+\sum_{l\in L_v}\frac{x-1}{x^{\mathfrak{d}(l)}-1}).
\end{align*}
\end{proof}

\begin{corollary}\label{zero forcing number of tree}
Let $T$ be a tree with $n$ vertices and $p$ many pendant vertices. Then $Z_C(T)=\lvert V({\rm st}(T))\rvert+p-\lvert S^T\rvert$. In addition, we have ${Z}_C(T;Z_C(T))=\underset{v\in S^T}{\prod}d_{{\rm st}(T)}(v)$.
\end{corollary}

\begin{proof}
From Theorem \ref{Thm: zero forcing polynomial of trees} it can be seen that $\frac{Z_C^T(x)}{x^{\rvert V({\rm st}(T))\lvert+p-\lvert S^T\rvert}}$ is a polynomial in $x$. Also, it is easy to see that $\frac{Z_C^T(x)}{x^{\rvert V({\rm st}(T))\lvert+p-\lvert S^T\rvert}}\Big\vert_{x=0}=\underset{v\in S^T}{\prod}\lvert L_v\rvert=\underset{v\in S^T}{\prod}d_{{\rm st}(T)}(v)$, which is nonzero. Hence, the smallest exponent of $x$ in $Z_C^{T}(x)$ with nonzero coefficient is $\lvert V({\rm st}(T))\rvert+p-\lvert S^T\rvert$, and the corresponding coefficient is $\underset{v\in S^T}{\prod}d_{{\rm st}(T)}(v)$. Thus, the result follows. 
\end{proof}

As a simple consequence of the above result, we have the following.
\begin{corollary}
Let $T_1$ and $T_2$ be two trees with $n$ vertices and $p$ many pendant vertices. Suppose both $T_1$ and $T_2$ have exactly two vertices $v_1,\ v_1'$ and $v_2,\ v_2'$ respectively, with degree more than $2$. Then ${\rm dist}(v_1,v_1')\geq{\rm dist}(v_2,v_2')$ implies $Z_C(T_1)\geq Z_C(T_2)$.
\end{corollary}

The following results are simple applications of results obtained in this section.
\begin{corollary}\label{computing connected zero forcing polynomial of special starlike}
Let $\mathcal{S}_{1,\ldots,1,n-p}$ be a starlike graph on $n$ vertices and $p$ many pendant vertices. Then $Z_C^{\mathcal{S}_{1,\ldots,1,n-p}}(x)=px^p+px^{p+1}+\cdots+px^{n-1}+x^n$.  
\end{corollary}

\begin{proof}
Consider the starlike graph $\mathcal{S}_{1,\ldots,1,n-p}$. Then from Lemma \ref{connected_zero forcing_starlike}, we have $$Z_C^{\mathcal{S}_{1,\ldots,1,n-p}}(x)=\frac{x^k(x^{n-p}-1)}{x-1}(x+p-1+\frac{x-1}{x^{n-p}-1})$$ $$=x^k(p\frac{x^{n-p}-1}{x-1}+x^{n-p})=x^k(p+px+\cdots+px^{n-p-1}+x^{n-p}).$$ Since degree of $Z_C^{\mathcal{S}_{1,\ldots,1,n-p}}(x)$ is $n$, we have $k=p$. Thus, the result holds.    
\end{proof}

\begin{corollary}\label{computing connected zero forcing polynomial of H}
 Let $2\leq k\leq p-2$ be a positive integer. Then $\mathcal{H}^n_{k,p-k}\in \mathfrak{T}_{n,p}$ and $Z_C^{\mathcal{H}^n_{k,p-k}}(x)=k(p-k)x^{n-2}+px^{n-1}+x^n$.  
 \end{corollary}
 \begin{proof}
    For $2\leq k\leq p-2$, consider the graph $\mathcal{H}^n_{k,p-k}$. Then from Theorem \ref{Thm: zero forcing polynomial of trees} we have $$Z_c^{\mathcal{H}^n_{k,p-k}}=x^k(x+k)(x+p-k)=x^{k+2}+px^{k+1}+k(p-k)x^k.$$ Since degree of $Z_c^{\mathcal{H}^n_{k,p-k}}$ is $n$, we have $k=n-2$. This concludes the proof.
\end{proof}

\section{Log-concavity of the polynomial}\label{Sec: Log}

\begin{lemma}[\cite{Br}; \cite{Ho}; \cite{St}, Prop.~2 ]\label{Lemma: Product of log-concave}
If $f,g$ are log-concave polynomials with nonnegative coefficients and no internal zeros, then so is $fg$.
\end{lemma}

From Proposition \ref{Prop:zero_forcing_polynomial_of_trees_structure}, it can be seen that the connected zero forcing polynomial of any finite tree is a product of restricted connected zero forcing polynomials of paths, and connected zero forcing polynomial of starlike graphs. Now, for paths, the restricted connected zero forcing polynomial of type $1$ has all coefficients $1$, hence it is log-concave.

\begin{lemma}\label{Lemma: log-concave path}
Let $T$ be a path of length at least $2$, and $u$ be an internal vertex of $T$. Then $Z_C^{T,u}(x)$ is log-concave with nonnegative coefficients and no internal zeros.
\end{lemma}
\begin{proof}
Let $T$ be a path of length at least $2$, with $u$ being an internal vertex. Let the two endpoints on the path be at a distance of $m$ and $n$ from $u$. Without loss of generality, we assume that $m\leq n$. Now, Lemma \ref{connected_zero_forcing_path} and some simple calculations yield that $Z_C^{T,u}(x)=x^2(x[m+1][n+1]+[m+1]+[n+1])$. Consider the polynomial $f(x)=x[m+1][n+1]+[m+1]+[n+1]$. Then it can be calculated that $f(x)=x(\sum_{i=0}^{m-1}(i+1)x^i+(m+1)\sum_{i=m}^nx^i+\sum_{i=1}^m(m+1-i)x^{n+i})+2\sum_{i=0}^mx^i+\sum_{i=m+1}^nx^i$, which can be further simplified to $f(x)=\sum_{i=0}^{m-1}(i+2)x^i+(m+2)\sum_{i=m}^nx^i+\sum_{i=1}^{m+1}(m+2-i)x^{n+i}$. Now, the coefficients of this polynomial form a sequence of the form $2,3,\ldots, m+1,m+2,\ldots,m+2,m+1,\ldots,3,2,1$, which is log-concave. Thus, the result holds.  
\end{proof}
Throughout the rest of the section, we will prove that the connected zero forcing polynomial of type $3$, i.e.~the for the starlike graph with more than two pendant vertices is log-concave. In order to prove that, let us first define some notations and present some important lemmas.

For a nonzero polynomial $f(x)=\sum_j f_jx^j$ with nonnegative coefficients, define $\lo(f):=\min\{j:\ f_j\neq 0\}$. So, $f=x^{\lo(f)}\hat f$ for some polynomial $\hat f$, with $\hat f(0)\ne 0$. For $N\in\mathbb Z$ we define the set
\[
  \U_N=\{0\}\cup\bigl\{f\neq 0~\vert ~f(x)=x^{\lo(f)}\hat f(x),~\hat f\text{ is symmetric and unimodal, and }\lo(f)+\deg (f)=N\bigr\},
\]
where $f$ ranges over polynomials with nonnegative coefficients. Equivalently, a nonzero polynomial $f(x)$ lies in $\U_N$ if and only if $f(x)=x^ag(x)$, for a positive integer $a$, and a symmetric unimodal polynomial $g(x)$ with nonnegative coefficients, such that $2a+\deg(g)=N$. For example, $1\in\U_0$ and $x^a[n]\in\U_{2a+n-1}$ for all integers $a\ge 0$, $n\ge 1$. Note that every nonzero $f\in\U_N$ has $\deg f=N-\lo(f)\leq N$.

\begin{lemma}\label{lem:U}
Let $M,N,m\in\mathbb Z$. Then we have the following;
\begin{enumerate}
\item[\rm (i)] $\U_N$ is closed under addition and under multiplication by nonnegative real numbers,
\item[\rm (ii)] If $f\in\U_M$ and $g\in\U_N$, then $fg\in\U_{M+N}$,
\item[\rm (iii)] If $f\in\U_{2m}$, then $f_{m-1}\leq f_m$.
\end{enumerate}
\end{lemma}

\begin{proof}
Let $f\in \U_N$. Write $f(x)=\sum_{i=0}^Nf_ix^i$, where $f_i=0$ for $0\leq i< \lo(f)$ and $\deg(f)<i\leq N$. Then we have $f_j=f_{N-j}$, for $0\leq j\leq N$. Also, since $f$ is unimodal, we have $f_0\leq f_1\leq \cdots\leq f_{\frac{N}{2}}\geq\cdots f_{N-1}\geq f_N$, if $N$ is even, and $f_0\leq f_1\leq \cdots\leq f_{\lfloor\frac{N}{2}\rfloor}=f_{\lfloor\frac{N}{2}\rfloor+1}\geq\cdots f_{N-1}\geq f_N$, if $N$ is odd. 
\begin{enumerate}
    \item[Proof of (i):] Consider $f,g\in \U_N$ is the above form. Then we have
    $f_j=f_{N-j}$ and $g_j=g_{N-j}$, for $0\leq j\leq N$. Thus, $f_j+g_j=f_{N-j}+g_{N-j}$ for $0\leq j\leq N$. So, $\frac{(f+g)(x)}{x^{\lo(f+g)}}$ is symmetric. Also, without loss of generality, assume that $N$ is even. Then we have $f_0\leq f_1\leq \cdots\leq f_{\frac{N}{2}}\geq\cdots f_{N-1}\geq f_N$ and $g_0\leq g_1\leq \cdots\leq g_{\frac{N}{2}}\geq\cdots g_{N-1}\geq g_N$. Thus, we have $f_0+g_0\leq f_1+g_1\leq \cdots\leq f_{\frac{N}{2}}+g_{\frac{N}{2}}\geq\cdots f_{N-1}+g_{N-1}\geq f_N+g_N$. So $f+g$ is unimodal. Also, the similar holds if $N$ is odd. 

    Using arguments similar to above, it can be shown that $\U_N$ is closed under multiplication by nonnegative real numbers.
    
    \item[Proof of (ii):] Let $f\in\U_M$ and $g\in\U_N$. So, $M=\deg(f)+\lo(f)$ and $N=\deg(g)+\lo(g)$. Also, we have $f(x)=x^{\lo(f)}\hat f(x)$ and $g(x)=x^{\lo(g)}\hat g(x)$, where $\hat{f}$ and $\hat{g}$ are symmetric and unimodal with nonnegative coefficients. Thus, from Lemma \ref{Lemma: Product symmetric unimodal} we have that $\hat{f}\cdot\hat{g}$ are symmetric and unimodal with nonnegative coefficients. Now, $f(x)g(x)=x^{\lo(f)+\lo(g)}\hat{f}(x)\hat{g}(x)$. Hence, $\lo(fg)+\deg(fg)=\lo(f)+\lo(g)+\deg(f)+\deg(g)=M+N$. So, $fg\in\U_{M+N}$.

    \item[Proof of (iii):] Let $N=2m$. Then a previous observation shows that $f_0\leq f_1\leq \cdots\leq f_{m-1}\leq f_m\geq f_{m+1}\geq\cdots f_{N-1}\geq f_N$. 
\end{enumerate}
This concludes the proof.
\end{proof}

Let $T=\Sk_{h_1,\dots,h_k}$ with center $v$ and legs of lengths $h_1,\dots,h_k\ge 1$. By the proof of Lemma \ref{connected_zero forcing_starlike}, a set $S\subseteq V(T)$ is a connected zero forcing set if and only if $S$ consists of $v$ together with the first $a_i$ vertices (counted from $v$) of the $i$-th leg, for some $0\leq a_i\leq h_i$, with $a_i=0$ for \emph{at most one} $i$. Define the set
\[
  J=J(h_1,\dots,h_k)=\bigl\{a\in\mathbb Z^k~\vert~ 0\le a_i\le h_i \text{ for }1\le i\le k,\ \#\{i:a_i=0\}\le 1\bigr\},
\]
and $J_m=\{a\in J:\ |a|=m\}$. Then
\begin{equation}\label{eq:model}
  \ZC^T(x)=x\sum_{a\in J}x^{|a|},\qquad\text{i.e.}\qquad \ZC^T(x)=x\sum_{i=o}^{\lvert h\rvert}\lvert J_i\rvert x^i.
\end{equation}

\subsection*{The fibre lemma}

For $s=(s_1,\dots,s_r)\in\mathbb Z_{\ge0}^r$, define
\[
  B(s):=\bigl\{b\in\mathbb Z^r~\vert~\ 0\leq b_i\leq s_i\text{ for }1\le i\le r,\ \#\{i:b_i=0\}\le1,\ \#\{i:b_i=s_i\}\le1\bigr\},
\]
and \(G_s(x)=\sum_{b\in B(s)}x^{|b|}\). Note that if $\#\{i:s_i=0\}\ge 2$, or $\#\{i:s_i=1\}\ge 3$, or $\#\{i:s_i=0\}= 1\text{ and }\#\{i:s_i=1\}\ge 1$, then $B(s)$ is empty. In that case, we consider $G_s(x)=0$. 

\begin{lemma}\label{lem:fibre}
For every $s\in\mathbb Z_{\ge0}^r$ we have $G_s\in\U_{|s|}$. In other words, either $G_s=0$, or
$G_s(x)=x^a\hat{G_s}(x)$ where $\hat{G_s}$ is a symmetric unimodal polynomial with nonnegative coefficients and
$2a+\deg(\hat{G_s})=|s|$.
\end{lemma}

\begin{proof}
Let $s=(s_1,\dots,s_r)\in\mathbb Z_{\ge0}^r$. If $B(s)=\emptyset$, then $G_s=0\in \U_{|s|}$. Otherwise, for $s_j\ge1$ define $M_j=x+x^2+\dots+x^{s_j-1}=x[s_j-1]$. If $s_j=1$ then $M_j=0$, and if $s_j\ge2$ then $M_j=x\cdot[s_j-1]$ with $[s_j-1]$ being symmetric unimodal of degree $s_j-2$, and $2\cdot1+(s_j-2)=s_j$. So $M_j\in\U_{s_j}$.

\emph{Case 1:} If some $s_i=0$, then $s_j\ge 2$ for all $i\neq j$, and $b_i=0=s_i$ for all $b\in B(s)$. So $b$ already has a coordinate $i$, which is equal to $0$ and equal to its maximum. Thus, all other coordinates must lie strictly inside their ranges. Hence, by Lemma~\ref{lem:U}(ii), we have $G_s=\prod_{j\ne i}M_j\in\U_{\sum_{j\ne i}s_j}=\U_{|s|}$.

\emph{Case 2:} If all $s_i\ge1$, then $\#\{i:s_i=1\}\le 2$. Relabel the coordinates of $s$ so that $s_1=\max_i s_i$. Now no coordinate can be simultaneously $0$ and maximal, so every $b\in B(s)$ has a well-defined \emph{type} $(\lambda,\mu)$, where $\lambda$ is the index $i$ with $b_i=0$ and $\mu$ is the index $i$ with $b_i=s_i$, $\lambda\ne\mu$. Also, if $b_i\neq 0$ for every $i$, we consider $\lambda=none$, and if $b_i\neq s_i$ for every $i$, we consider $\mu=none$. The elements of a given type have generating function $\prod_i\phi_i$ with $\phi_\lambda=1$,
$\phi_\mu=x^{s_\mu}$ and $\phi_i=M_i$ otherwise. Grouping the types we get
$G_s=\Gamma_1+\sum_{j\ge2}\Gamma_j+\sum_{2\leq j<j'}\Gamma_{jj'}$, where:
\begin{itemize}
\item $\Gamma_1$ collects the types (none,\,none), $(1,\text{none})$, $(\text{none},1)$:
  \[
    \Gamma_1=(M_1+1+x^{s_1})\prod_{l\ne1}M_l=[s_1+1]\prod_{l\ne1}M_l .
  \]
\item For $j\ge2$, $\Gamma_j$ collects $(j,\text{none})$, $(\text{none},j)$, $(1,j)$, $(j,1)$:
  \[
    \Gamma_j=\Bigl[(1+x^{s_j})M_1+x^{s_j}+x^{s_1}\Bigr]\prod_{l\notin\{1,j\}}M_l
            =(x+x^{s_j})\,[s_1]\prod_{l\notin\{1,j\}}M_l,
  \]
  using $M_1+x^{s_1}=x[s_1]$ and $M_1+1=[s_1]$.
\item For $2\leq j<j'$, $\Gamma_{jj'}$ collects $(j,j')$ and $(j',j)$:
  \[
    \Gamma_{jj'}=(x^{s_j}+x^{s_{j'}})\,M_1\prod_{l\notin\{1,j,j'\}}M_l .
  \]
\end{itemize}
Recall that $M_l\in\U_{s_l}$ for $1\le l\le r$. We will show that each $\Gamma$ lies in $\U_{|s|}$. 

\vspace{.3cm}
\noindent\underline{$\Gamma_1$:} $[s_1+1]$ is symmetric unimodal of degree $s_1$, so
$[s_1+1]\in\U_{s_1}$. Thus $\Gamma_1\in\U_{|s|}$ by Lemma~\ref{lem:U}(ii).

\vspace{.3cm}
\noindent\underline{$\Gamma_j$:} We have $(x+x^{s_j})[s_1]=x\cdot(1+x^{s_j-1})[s_1]$. Since
$0\leq s_j-1\leq s_1-1$, the exponent ranges of $[s_1]$ and $x^{s_j-1}[s_1]$ overlap. Now, a simple calculation shows that the coefficient sequence
of $(1+x^{s_j-1})[s_1]$ is
\[
  \underbrace{1,\dots,1}_{s_j-1},\ \underbrace{2,\dots,2}_{s_1-s_j+1},\ \underbrace{1,\dots,1}_{s_j-1}.
\]
This polynomial is symmetric unimodal of degree $s_1+s_j-2$, so
$(x+x^{s_j})[s_1]\in\U_{2+s_1+s_j-2}=\U_{s_1+s_j}$, and $\Gamma_j\in\U_{|s|}$ by
Lemma~\ref{lem:U}(ii).

\vspace{.3cm}
\noindent\underline{$\Gamma_{jj'}$:} If $s_1=1$ then $M_1=0$ and $\Gamma_{jj'}=0$. Let $s_1\ge2$.
Without loss of generality, we assume $s_j\leq s_{j'}$, and denote
$\delta=s_{j'}-s_j$. Then
\[
  (x^{s_j}+x^{s_{j'}})M_1=x^{s_j+1}(1+x^{\delta})[s_1-1],
\]
and $0\le\delta\leq s_1-1$ because $1\leq s_j\leq s_{j'}\leq s_1$. So the exponent ranges
of $[s_1-1]$ and $x^{\delta}[s_1-1]$ overlap or are adjacent. Now, a simple calculation shows that the coefficient sequence of $(1+x^{\delta})[s_1-1]$ is
\[
  \underbrace{1,\dots,1}_{\delta},\ \underbrace{2,\dots,2}_{s_1-1-\delta},\ \underbrace{1,\dots,1}_{\delta}.
\]
This polynomial is symmetric unimodal of degree $\delta+s_1-2$, so
$(x^{s_j}+x^{s_{j'}})M_1\in\U_{2(s_j+1)+\delta+s_1-2}=\U_{s_1+s_j+s_{j'}}$, and
$\Gamma_{jj'}\in\U_{|s|}$ by Lemma~\ref{lem:U}(ii).

\smallskip
Finally, since $G_s=\Gamma_1+\sum_{j\ge2}\Gamma_j+\sum_{2\leq j<j'}\Gamma_{jj'}$, by Lemma~\ref{lem:U}(i) we have $G_s\in\U_{|s|}$.
\end{proof}

\subsection*{Log-concavity for starlike graphs}

\begin{proposition}\label{thm:main}
Let $k\ge 3$ and $h_1,\dots,h_k\ge1$. Then for the starlike graph $T=\Sk_{h_1,\dots,h_k}$, the connected zero forcing polynomial $\ZC^T(x)$ is log-concave with no internal zeros.
\end{proposition}

\begin{proof}
Consider the starlike graph $\Sk_{h_1,\dots,h_k}$. By \eqref{eq:model} it suffices to show that the coefficient sequence $(|J_m|)_m$ is log-concave with no internal zeros.

\underline{No internal zeros:} Since $k\ge 3$, an element of $J$ has at least $k-1$ positive coordinates, so $J_m=\emptyset$ for $m\le k-2$. Now, starting from $(0,1,\dots,1)\in J_{k-1}$ and raising coordinates one
unit at a time up to $(h_1,\dots,h_k)$, we stay inside $J$ and visit every level in between. Thus, $(|J_m|)_m$ has no internal zeros.

\underline{Log-concavity:} Fix an $m$ so that $k\leq m\leq \lvert h\rvert-1$. For $\sigma\in\mathbb Z^k$ define the \emph{fibre}
\[
  F(\sigma)=\{a\in J:\ \sigma-a\in J\},\qquad c_\sigma(t)=\#\{a\in F(\sigma):\ |a|=t\}.
\]
Now, $|\sigma|=2m$, and any $a\in F(\sigma)$ with $|a|=m$ (resp.\ $m-1$) implies that
$b=\sigma-a$ automatically lies in $J_m$ (resp.\ $J_{m+1}$). So, counting ordered pairs $(a,b)$ according to $\sigma=a+b$ gives
\begin{equation}\label{eq:pairs}
  |J_m|^2=\sum_{|\sigma|=2m}c_\sigma(m),\qquad |J_{m-1}|\,|J_{m+1}|=\sum_{|\sigma|=2m}c_\sigma(m-1),
\end{equation}
Thus, it suffices to prove
$c_\sigma(m)\ge c_\sigma(m-1)$ for every $\sigma$ with $|\sigma|=2m$, and we may assume
$0\le\sigma_i\le2h_i$ for all $i$ (otherwise $F(\sigma)=\emptyset$).

Now $a\in F(\sigma)$ means: $\ell_i\leq a_i\leq u_i$ with $\ell_i=\max(0,\sigma_i-h_i)$ and
$u_i=\min(h_i,\sigma_i)$; $a$ has at most one zero coordinate; and $\sigma-a$ has at most one zero
coordinate, i.e.~$a_i=\sigma_i$ for at most one $i$. Let $I=\{i:\ \sigma_i\leq h_i\}$. For $i\in I$
the range of $a_i$ is $[0,\sigma_i]$, and $\#(i:a_i=0)\le 1$ and $\#(i:a_i=\sigma_i)\le 1$. For $i\notin I$ we have $\ell_i=\sigma_i-h_i\ge1$ and $u_i=h_i<\sigma_i$, so
$a_i\ne0$ and $a_i\ne\sigma_i$. Consequently, we have
\[
  \sum_t c_\sigma(t)\,x^t=\sum_{a\in F(\sigma)}x^{|a|}
  =\Bigl(\prod_{i\notin I}x^{\ell_i}[u_i-\ell_i+1]\Bigr)\cdot G_{(\sigma_i)_{i\in I}}(x).
\]
For $i\notin I$ we have $u_i-\ell_i+1\ge1$ (as $\sigma_i\le2h_i$), so $[u_i-\ell_i+1]$ is
symmetric unimodal of degree $u_i-\ell_i$. Since $2\ell_i+(u_i-\ell_i)=\ell_i+u_i=\sigma_i$, the
factor $x^{\ell_i}[u_i-\ell_i+1]$ lies in $\U_{\sigma_i}$. By Lemma~\ref{lem:fibre},
$G_{(\sigma_i)_{i\in I}}\in\U_{\sum_{i\in I}\sigma_i}$. So, applying Lemma~\ref{lem:U}(ii), we have
\[
  \sum_t c_\sigma(t)\,x^t=\Bigl(\prod_{i\notin I}x^{\ell_i}[u_i-\ell_i+1]\Bigr)\cdot G_{(\sigma_i)_{i\in I}}(x)\in\U_{|\sigma|}=\U_{2m},
\]
and Lemma~\ref{lem:U}(iii) gives $c_\sigma(m-1)\leq c_\sigma(m)$. Summing over $\sigma$ in
\eqref{eq:pairs} yields $|J_m|^2\ge|J_{m-1}|\,|J_{m+1}|$. This concludes the proof.
\end{proof}

So, throughout this section we have shown that for a finite tree $T$ both the polynomials $Z_C^{T_v,v}(x)$ for $v\in S_1^T\cup S_2^T$ and $Z_C^{T_w}(x)$ for $w\in S_3^T$ are log concave with no internal zeroes. Thus, with Lemma \ref{Lemma: Product of log-concave} and Proposition \ref{Prop:zero_forcing_polynomial_of_trees_structure} in hand, we conclude the following.

\begin{theorem}\label{Thm: Log concave trees}
For any finite tree $T$, the connected zero forcing polynomial $Z_C^T(x)$ is log-concave with no internal zeroes.    
\end{theorem}
As a corollary of the above result, we resolve an open question posed by Newton and Brimkov in \cite{Newton},
\begin{corollary}\label{Cor: Unimodality}
For any finite tree $T$, the connected zero forcing polynomial $Z_C^T(x)$ is unimodal.   
\end{corollary}

\section{Bounds for the connected zero forcing polynomial of trees}
In \cite{boyer}, the authors conjectured that, if we consider the path $p_n$ with $n$ vertices, then for any graph $G$ with $n$ vertices, $Z(p_n;i) \geq Z(G,i)$ for all $i\geq 1$. In \cite{ka} and \cite{german}, the authors included  some classes of graphs in support of this conjecture of \cite{boyer}. But for connected zero forcing sets the scenario is totally different. For example, consider the following graphs $\mathcal{S}_{1,1,1}$ and $P_4$. Then

$Z_C^{\mathcal{S}_{1,1,1}}(x)=x^4+3x^3+0.x^2+0.x$ and $Z_C^{P_4}(x)=x^4+2x^3+3.x^2+2x$.\\
Clearly, $Z_C(P_4,3)\ngeq Z_C(\mathcal{S}_{1,1,1},3)$. In fact, the polynomials $Z_C^{\mathcal{S}_{1,1,1}}(x)$ and $Z_C^{P_4}(x)$ are not comparable. Though the polynomial $Z(\mathcal{S}_{1,1,\cdots,1},i)\leq Z(G,i)$ for $1\leq i\leq n$ and for any graph $G$ with $n$ vertices. But $Z_C^{\mathcal{S}_{1,1,\cdots,1}}(x)$ and $Z_C^{P_n}(x)$ are even not comparable. In fact, if we add a new vertex $v$ by an edge with some vertex $u$ of some tree $T$, then we have the following Lemma \ref{add.ver}.
\begin{lemma}\label{add.ver}
Let $T$ be a tree which is not a path, and $u\in V(T)$. Let $T^u\cup \{v\}$ denote the new tree formed after attaching a new pendant vertex $v$ with $u$ by an edge. Then we have the following:
\begin{itemize}
    \item If $u\in P(T)$, then ${Z}_C(T^u\cup\{v\};i)={Z}_C(T;i)+{Z}^{u}_C(T;i-1)$, for $i> 2$,
    \item If $u\notin P(T)$, then ${Z}_C(T^u\cup\{v\};i)={Z}^{N_T[u]}_C(T;i)+{Z}^{u}_C(T;i-1)$.
\end{itemize}
\end{lemma}
\begin{proof}
    $T$ is a tree with $n-1$ vertices and $p-1$ many pendant vertices, and $u\in V(T)$. $T^u\cup \{v\}$ is the new tree formed after attaching a new pendant vertex $v$ with $u$ by an edge. It is easy to see that any connected zero forcing set of $T^u\cup \{v\}$ containing $v$ also contains $u$, since $T$ is not a path. In addition, any connected zero forcing set of $T^u\cup \{v\}$ which does not contain $v$ is a connected zero forcing set of $T$.
    \begin{itemize}
        \item[\underline{Case 1:}] If $u\in P(T)$, then any connected zero forcing set of $T$ will also be a connected zero forcing set of $T^u\cup \{v\}$. Thus, we have ${Z}_C(T^u\cup\{v\};i)={Z}_C(T;i)+{Z}^{u}_C(T;i-1)$, for $i> 2$. 

        \item[\underline{Case 2:}] If $u\notin P(T)$, for any connected zero forcing set of $T^u\cup \{v\}$ which does not contain $v$, the only way $v$ can be colored is when $v$ gets colored by $u$, and $u$ can color $v$ only when all the vertices in $N_T[u]$ are colored. Since $T^u\cup \{v\}$ is a tree, a connected zero forcing set of $T^u\cup \{v\}$ can induce coloring on all the vertices in $N_T[u]$ only when the connected zero forcing set itself contains the set $N_T[u]$. On the other hand, it is easy to see that all the connected zero forcing set of $T$ which contains $N_T[u]$ are also connected zero forcing sets of $T^u\cup \{v\}$. Thus, we have ${Z}_C(T^u\cup\{v\};i)={Z}^{N_T[u]}_C(T;i)+{Z}^{u}_C(T;i-1)$. 
    \end{itemize}
\end{proof}
From the above discussion, it seems that there is no relation between the coefficients of the connected zero forcing polynomials of the trees. But that is not the scenario if we fix the number of pendant vertices of the trees. In fact, when the number of pendant vertices $p$ is fix, then for trees in $\mathfrak{T}_{n,p}$, we prove interesting results in Sections \ref{lower.b} and \ref{upper.b}.

\subsection{Lower bounds for the connected zero forcing polynomial of trees}\label{lower.b}

In this section, we establish several important and interesting results concerning the connected zero forcing polynomials of trees, namely, Proposition \ref{lower bound for any tree} and Proposition \ref{the case for nontrivial st}. In particular, we determine lower bounds for the connected zero forcing polynomials of trees belonging to $\mathfrak{T}_{n,p}$. In these two propositions, we show that although the connected zero forcing polynomials of the trees $\mathcal{S}_{1,1,\ldots,1,n-p}$ and $\mathcal{H}^{n}_{2,p-2}$ in $\mathfrak{T}_{n,p}$ are incomparable, every other tree $T\in\mathfrak{T}_{n,p}$ satisfies
$$
Z_C^T(x)\geq Z_C^{\mathcal{S}_{1,1,\ldots,1,n-p}}(x)
$$
and
$$
Z_C^T(x)\geq Z_C^{\mathcal{H}^{n}_{2,p-2}}(x).
$$
Thus, the two trees $\mathcal{S}_{1,1,\ldots,1,n-p}$ and $\mathcal{H}^{n}_{2,p-2}$ are precisely the two minimal elements of $\mathfrak{T}_{n,p}$ with respect to the ordering induced by their connected zero forcing polynomials.
\begin{lemma}\label{some coeffecients}
Let $T\in \mathfrak{T}_{n,p}$. Then we have the following:
\begin{enumerate}
    \item ${Z}_C(T;n-1)=p$,
    \item If $\lvert S^T\rvert\geq 2$, then there exists an integer $m$ with $2\leq m\leq p-2$, such that ${Z}_C(T;n-2)\geq m(p-m)$. 
\end{enumerate}
\end{lemma}

\begin{proof}
Consider the reciprocal of the polynomial $Z_C^{T}(x)$ to be $Z_C^{T*}(x)$. Then theorem \ref{Thm: zero forcing polynomial of trees} and some simple calculations yield $$Z_C^{T*}(x)=\underset{l\in P(T)}{\prod}(\frac{x^{\mathfrak{d}(l)}-1}{x-1})\times\underset{v\in S^T}{\prod}(1+\underset{l\in L_v}{\sum}(x^{\mathfrak{d}(l)}\frac{x-1}{x^{\mathfrak{d}(l)}-1})).$$ So, ${Z}_C(T;n-1)$ is the coefficient of $x$ in $Z_C^{T*}(x)$, and ${Z}_C(T;n-2)$ is the coefficient of $x^2$ in $Z_C^{T*}(x)$.
\begin{enumerate}
    \item It can be calculated that $\frac{d}{dx}Z_C^{T*}(x)\vert_{x=o}=\underset{\substack{l\in P(T),\\\mathfrak{d}(l)=1}}{\sum}1+\underset{\substack{l\in P(T),\\\mathfrak{d}(l)\geq 2}}{\sum}1=p$. So, the coefficient of $x$ in $Z_C^{T*}(x)$, i.e.~${Z}_C(T;n-1)=p$. 
    \item If $\lvert S^T\rvert\geq 2$, it is easy to see that $\lvert S_2^T\rvert\geq 2$. Choose $u\in S_2^T$. Clearly $2\leq \lvert L_u\rvert\leq p-2$. Now, $$Z_C^{T*}(x)=(\underset{l\in L_u}{\prod}(\frac{x^{\mathfrak{d}(l)}-1}{x-1})\times(1+\underset{l\in L_u}{\sum}(x^{\mathfrak{d}(l)}\frac{x-1}{x^{\mathfrak{d}(l)}-1})))\times$$ $$(\underset{l\in P(T)\setminus L_u}{\prod}(\frac{x^{\mathfrak{d}(l)}-1}{x-1})\times\underset{v\in S^T\setminus\{u\}}{\prod}(1+\underset{l\in L_v}{\sum}(x^{\mathfrak{d}(l)}\frac{x-1}{x^{\mathfrak{d}(l)}-1}))).$$ Then the previous case yields that the coefficient of $x$ in $\underset{l\in L_u}{\prod}(\frac{x^{\mathfrak{d}(l)}-1}{x-1})\times(1+\underset{l\in L_u}{\sum}(x^{\mathfrak{d}(l)}\frac{x-1}{x^{\mathfrak{d}(l)}-1}))$ is $\lvert L_u\rvert$ and the coefficient of $x$ in $\underset{l\in P(T)\setminus L_u}{\prod}(\frac{x^{\mathfrak{d}(l)}-1}{x-1})\times\underset{v\in S^T\setminus\{u\}}{\prod}(1+\underset{l\in L_v}{\sum}(x^{\mathfrak{d}(l)}\frac{x-1}{x^{\mathfrak{d}(l)}-1}))$ is $p-\lvert L_u\rvert$. Hence, we have the coefficient of $x^2$ in $Z_C^{T*}(x)$ is greater than or equal to $\lvert L_u\rvert(p-\lvert L_u\rvert)$, where $2\leq \lvert L_u\rvert\leq p-2$. This concludes the proof.
\end{enumerate}  
\end{proof}

The following observation is useful.
\begin{observation}\label{Obs: polynomial coefficient observation}
Let $H=\{h_1,h_2,\ldots,h_k\}$ be a multiset of $k$ positive integers. Consider the polynomial $g_H(x)=\underset{i=1}{\overset{k}{\prod}}\underset{j=0}{\overset{h_i}{\sum}}x^j$. Then the coefficient of $x^l$ in $g(x)$ is bounded below by $k$, for $1\leq l\leq\underset{i=1}{\overset{k}{\sum}}h_i-1$.    
\end{observation}
\begin{proof}
 The proof proceeds by induction on $k$. Clearly, the result is true for $k=1$. Let the result be true for $k=m-1$ for a positive integer $m\geq 2$. Consider the multiset $H=\{h_1,h_2,\ldots,h_m\}$ of $m$ positive integers. Define $h=h_1+h_2+\cdots h_{m-1}$. Then $$g_H(x)=(\prod_{i=1}^{m-1}\sum_{j=0}^{h_i}x^j)(\sum_{j=0}^{h_m}x^j).$$ Let $\prod_{i=1}^{m-1}\sum_{j=0}^{h_i}x^j=\sum_{i=0}^hc_ix^i$, where $c_0=1=c_{h_1+h_2+\cdots h_{m-1}}$ and $c_i\geq m-1$ for $1\leq i\leq h-1$. Then coefficient of $x^l$ in $g_H(x)$ is $\underset{\substack{i+j=l,\\0
 \leq i\leq h,\\ 0\leq j\leq h_k}}{\sum}c_i$. Now, for $1\leq l\leq h+h_k-1$ we can choose integers $i_1,j_1$ such that $0
 \leq i_1<i_1+1\leq h$, $0\leq j_1-1<j_1\leq h_k$ and $i_1+j_1=l$. So, $\underset{\substack{i+j=l,\\0
 \leq i\leq h,\\ 0\leq j\leq h_k}}{\sum}c_i\geq c_{i_1}+c_{i_1+1}\geq m$. Hence the statement holds for $k=m$. This concludes the proof. 
\end{proof}

We will now find a lower bound for the coefficients of the connected zero forcing polynomial of any $T\in \mathfrak{T}_{n,p}$. Before going into the main result, we will first explore such lower bounds for $Z_C^{T_v,v}(x)$, when $v\in S_1^T\cup S_2^T$, and for $Z_C^{T_v}(x)$, when $v\in S_3^T$. Note that for $v\in S_1^T$, $Z_C^{T_v,v}(x)$ is of the form $x+x^2+\cdots+x^{k+1}$, where $k$ is the length of the path $T_v$. Hence, each nonzero coefficient is bounded below by the number of pendant vertices of $T$, which $T_v$ contains. The next result shows that such phenomena also holds for $v\in S_2^T$ or and $v\in S_3^T$.

\begin{lemma}\label{Lemma: Lower bound for starlike}
Let $T$ be a tree on $n$ vertices, such that either $T$ is a path of length more than $1$ with an internal vertex $v$, or a starlike graph. Then we have the following,
\begin{itemize}
    \item If $T$ is a path, then each nonzero coefficient of $Z_C^{T,v}(x)$ other than the first and the last one, is bounded below by $2$.
    \item If $T$ is a starlike graph, then $Z_C(T;i)\geq p(T)$, for $p(T)< i< n$.
\end{itemize}
\end{lemma}
\begin{proof}
Let $T$ be either a path with an internal vertex $v$ or a starlike graph. Further, if $T$ is a starlike graph, we denote the central vertex of $T$ by $v$. In either cases, $\mathfrak{d}(l)$ will denote the distance between $l$ and $v$, for an $l\in P(T)$. Then, Lemma \ref{connected_zero_forcing_path} implies that $Z_C^{T,v}(x)=x^2\frac{\prod_{l\in P(T)}(x^{\mathfrak{d}(l)}-1)}{(x-1)^2}(x+\underset{l\in P(T)}{\sum}\frac{x-1}{x^{\mathfrak{d}(l)}-1})$, when $T$ is a path, and Lemma \ref{connected_zero forcing_starlike} implies that $Z_C^T(x)=x^{p(T)}(\underset{l\in P(T)}{\prod}\frac{x^{\mathfrak{d}(l)}-1}{x-1})(x+\underset{l\in P(T)}{\sum}\frac{x-1}{x^{\mathfrak{d}(l)}-1})$, when $T$ is starlike. Now, a few steps of simple calculations yield that in both the cases, the polynomial is of the form $x^{p(T)}(x\underset{l\in P(T)}{\prod}\underset{i=0}{\overset{\mathfrak{d}(l)-1}{\sum}}x^i+\underset{l\in P(T)}{\sum}(\underset{\substack{l'\in P(T),\\ l'\neq l}}{\prod}\underset{i=0}{\overset{\mathfrak{d}(l')-1}{\sum}}x^i))$.  Consider the polynomial \[f(x):=x\underset{l\in P(T)}{\prod}\underset{i=0}{\overset{\mathfrak{d}(l)-1}{\sum}}x^i+\underset{l\in P(T)}{\sum}(\underset{\substack{l'\in P(T),\\ l'\neq l}}{\prod}\underset{i=0}{\overset{\mathfrak{d}(l')-1}{\sum}}x^i).\] We will show that each nonzero coefficient of $f(x)$ is bounded below by $p(T).$ 

Define $P_1(T):=\{l\in P(T)~\vert~\mathfrak{d}(l)=1\}$, and $P_2(T):=\{l\in P(T)~\vert~\mathfrak{d}(l)\geq 2\}$. So, $P(T)=P_1(T)\sqcup P_2(T)$. Then \[f(x)=(x+\lvert P_1(T)\rvert)(\underset{l\in P_2(T)}{\prod}\underset{i=0}{\overset{\mathfrak{d}(l)-1}{\sum}}x^i)+\underset{l\in P_2(T)}{\sum}(\underset{\substack{l'\in P_2(T),\\ l'\neq l}}{\prod}\underset{i=0}{\overset{\mathfrak{d}(l')-1}{\sum}}x^i).\] 
Let the polynomial $\underset{l\in P_2(T)}{\prod}\underset{i=0}{\overset{\mathfrak{d}(l)-1}{\sum}}x^i$ be of the form $\underset{j=0}{\overset{k}{\sum}}c_ix^i$, where $k=\underset{l\in P_2(T)}{\sum}(\mathfrak{d}(l)-1)$. Then \[(x+\lvert P_1(T)\rvert)(\underset{l\in P_2(T)}{\prod}\underset{i=0}{\overset{\mathfrak{d}(l)-1}{\sum}}x^i)=\lvert P_1(T)\rvert+\sum_{i=1}^k(\lvert P_1(T)\rvert c_i+c_{i-1})x^i+x^{k+1}.\]
Clearly $\mathfrak{d}(l)-1$ are all positive integers for all $l\in P_2(T)$. So we have $c_0=1=c_k$, and by Observation \ref{Obs: polynomial coefficient observation} we have $c_i\geq \lvert P_2(T)\rvert$, for $1\leq i\leq k-1$. Hence, $\lvert P_1(T)\rvert c_i+c_{i-1}\geq \lvert P_1(T)\rvert+\lvert P_2(T)\rvert=\lvert P(T)\rvert$, for $2\leq i\leq k$. If $\lvert P_1(T)\rvert\geq 1$, then $\lvert P_1(T)\rvert c_1+c_0\geq \lvert P_1(T)\rvert+\lvert P_2(T)\rvert=\lvert P(T)\rvert$, since $c_0=1$ and $c_1\geq \lvert P_2(T)\rvert$. If $\lvert P_1(T)\rvert=0$, then $P_2(T)=P(T)$. Hence, the coefficient of $x$ in $\underset{l\in P(T)}{\sum}(\underset{\substack{l'\in P(T),\\ l'\neq l}}{\prod}(\underset{i=0}{\overset{\mathfrak{d}(l')-1}{\sum}}x^i))$ is $\lvert P(T)\rvert(\lvert P(T)\rvert-1)\geq \lvert P(T)\rvert$. Hence, the coefficient of $x^j$ in $f(x)$ is greater than or equal to $p(T)$, for $1\leq j\leq k$. 
\end{proof}

Also, the following holds.
\begin{lemma}\label{Lemma: product lower bound}
 Let $f_1(x)=\sum_{i=0}^mc_ix^i$ and $f_2(x)=\sum_{j=0}^nd_jx^j$, $m,n\geq 2$ be two polynomials with positive integral coefficients, such that $c_i\geq a$ and $d_j\geq $, for $1\leq i\leq m-1$, $1\leq j\leq n-1$. If $f_1(x)f_2(x)=\sum_{i=0}^{m+n}l_ix^i$, then $l_i\geq a+b$ for $1\leq i\leq m+n-1$.
\end{lemma}
\begin{proof}
 $f_1(x)=\sum_{i=0}^mc_ix^i$ and $f_2(x)=\sum_{j=0}^nd_jx^j$ are two polynomials with positive integral coefficients, such that $c_i\geq a$ and $d_j\geq b$, for $1\leq i\leq m-1$, $1\leq j\leq n-1$. Since $f_1(x)f_2(x)=\sum_{i=0}^{m+n}l_ix^i$, we have $l_i=\underset{\substack{h+k=i \\ 0\leq h\leq m \\ o\leq k\leq n}}{\sum}c_hd_k$, for $0\leq i\leq m+n$. Now, for $1\leq l\leq m+n-1$ we can choose integers $h',k'$ such that $0
\leq h'<h'+1\leq m$, $0\leq k'-1<k'\leq n$ and $h'+k'=l$. So, we have $l_i=\underset{\substack{h+k=i \\ 0\leq h\leq m \\ o\leq k\leq n}}{\sum}c_hd_k\geq c_{h'}d_{k'}+c_{h'+1}d_{k'-1}\geq a+b$. Thus, the result holds.
\end{proof}

Note that the above result still holds if one or both of the polynomials are of the form $1+x$. So, as a consequence of the above two results, we now arrive at the following lower bound for the coefficients of the connected zero forcing polynomial of trees.

\begin{proposition}\label{lower bound for any tree}
Let $T\in \mathfrak{T}_{n,p}$ with $Z_C(T)=k$. Then ${Z}_C(T;i)\geq{Z}_C(\mathcal{S}_{1,1,\ldots,n-p};i)=p$ for $k<i<n$.   
\end{proposition}

\begin{proposition}\label{the case for nontrivial st}
Let $n,p\in\mathbb{N}$ with $p\geq 4$, $n\geq p+2$, and $T\in \mathfrak{T}_{n,p}$ such that $T\neq \mathcal{S}_{1,\ldots,1,n-p}$. Then ${Z}_C(T;i)\geq{Z}_C(\mathcal{H}^n_{2,p-2};i)$ for all $i\in \{1,2,\ldots,n\}$.   
\end{proposition}
\begin{proof}
Let $n,p\in\mathbb{N}$ with $p\geq 4$ and $n\geq p+2$. Let $T\in \mathfrak{T}_{n,p}$ such that $T\neq \mathcal{S}_{1,\ldots,1,n-p}$. When $S^T=\{u\}$ we have $$Z_C^{T*}(x)=\underset{l\in P(T)}{\prod}(\frac{x^{\mathfrak{d}(l)}-1}{x-1})+\underset{l\in L_v}{\sum}(x^{\mathfrak{d}(l)}\underset{l'\in P(T)\setminus\{l\}}{\prod}(\frac{x^{\mathfrak{d}(l)}-1}{x-1})).$$
Let $P(T)_1:=\{l\in P(T)~\vert~\mathfrak{d}(l)=1\}$, $P(T)_2:=\{l\in P(T)~\vert~\mathfrak{d}(l)=2\}$, $P(T)_3:=\{l\in P(T)~\vert~\mathfrak{d}(l)\geq 3\}$, and $\lvert P(T)_2\rvert+\lvert P(T)_3\rvert=t$. Clearly, $\lvert P(T)_1\rvert+t=p$ and $t\geq 2$ as $T\neq \mathcal{S}_{1,\ldots,1,n-p}$. Then $$Z_C^{T*}(x)=(1+\lvert P(T)_1\rvert x)\underset{l\in P(T)_2\cup P(T)_3}{\prod}(\frac{x^{\mathfrak{d}(l)}-1}{x-1})+x^2\underset{l\in P(T)_2}{\sum}(\underset{\substack{l'\in P(T)_2\cup P(T)_3, l'\neq l}}{\prod}(\frac{x^{\mathfrak{d}(l)}-1}{x-1}))$$ $$+\underset{l\in P(T)_3}{\sum}(x^{\mathfrak{d}(l)}\underset{\substack{l'\in P(T)_2\cup P(T)_3,\\ l'\neq l}}{\prod}(\frac{x^{\mathfrak{d}(l)}-1}{x-1})).$$
It is easy to calculate that the coefficient of $x^2$ in $Z_C^{T*}(x)$ is $t(\lvert P(T)_1\rvert+1)+\binom{t}{2}$. If $t=p$, then $$t(\lvert P(T)_1\rvert+1)+\binom{t}{2}=p+\binom{p}{2}\geq 2(p-2),$$
and if $t=p-1$, then $$t(\lvert P(T)_1\rvert+1)+\binom{t}{2}=p-1+\binom{p-1}{2}\geq 2(p-2),\ \mbox{since}\ p\geq 4.$$ 
If $2\leq t\leq p-2$, then $$t(\lvert P(T)_1\rvert+1)+\binom{t}{2}=t(p-\frac{t-1}{2})\geq t(p-t).$$
Now, a few simple calculations yield that $t(p-t)\geq 2(p-2)$, for $2\leq t\leq p-2$. Thus, from corollary \ref{some coeffecients} and \ref{computing connected zero forcing polynomial of H}, it can be seen that if $T\in \mathfrak{T}_{n,p}$ such that $T\neq \mathcal{S}_{1,\ldots,1,n-p}$, then ${Z}_C(T;i)\geq{Z}_C(\mathcal{H}^n_{2,p-2};i)$, for all $i\in \{1,2,\ldots,n\}$. This concludes the proof.
\end{proof}


\subsection{Upper bounds for the connected zero forcing polynomial of trees}\label{upper.b}
We begin this section with a definition and one of the most important results of this paper, namely, Theorem \ref{thm.uniform}. We then prove that, for $p=2$ and $3$, there always exists a unique maximal tree in $\mathfrak{T}_{n,p}$. However, for $p\geq 5$, we present several interesting examples demonstrating that, as the number of pendant vertices increases, there exists a certain threshold value of $n$ beyond which no tree $T\in\mathfrak{T}_{n,p}$ satisfies
$
Z_C(T',i)\leq Z_C(T,i)
$
for every $1\leq i\leq n$ and every $T'\in\mathfrak{T}_{n,p}$. In other words, for such values of $n$ and $p$, there is no tree in $\mathfrak{T}_{n,p}$ whose connected zero forcing polynomial coefficients simultaneously dominate those of every other tree in the family.
\begin{definition}
Let $T\in \mathfrak{T}_{n,p}$ and $p_1,\ p_2$ be any two pendant vertices of $T$, such that the distance of $p_1$ and $p_2$ from their seeding vertices are $d_1$ and $d_2$ respectively (here the seeding vertices of $p_1$ and $p_2$ can be same). Then $T$ is said to be an uniform tree if  $|d_1-d_2|$ is either $0$ or $1$ for any two pendant vertices $p_1$ and $p_2$ in $T$.
   
\end{definition}

\begin{theorem}\label{thm.uniform}
Let $T'$ be any tree in $\mathfrak{T}_{n,p}$. Also let $T'$ has the seeding vertices $\{s_1,\cdots, s_k\}$ and for $1\leq j \leq k$, each of the seeding vertex $s_j$ has $a_j$ many seeds. Then there always exists a uniform tree $T$ such that $Z_C(T',i)\leq Z_C(T,i)$ for $1\leq i \leq n$.
\end{theorem}
\begin{proof}
    Let us consider the polynomials $(1+x+x^2+\cdots x^m)$ and $(1+x+x^2+\cdots x^r)$, where $m\leq r$. Let
    $$\sum_{i=0}^mx^i \sum_{j=0}^rx^j=\sum_{k=0}^{m+r}c_tx^t.$$
    Then we know that the coefficient of $x^t$ in the product of these two polynomials is given by 
\begin{equation}\label{equ.poly}
   c_t=
\begin{cases}
t+1, & \text{if } 0 \leq t \leq m,\\[4pt]
m+1, & \text{if } m \leq t \leq r,\\[4pt]
m+r-t+1, & \text{if } r \leq t \leq m+r,\\[4pt]
0, & \text{otherwise}.
\end{cases} 
\end{equation}
From Theorem \ref{Thm: zero forcing polynomial of trees}, we know that if $T\in \mathfrak{T}_{n,p}$ and $Z_C(T)=k$. Then 
\begin{equation}\label{equ.zerof}
    Z_C^{T}(x)=\frac{x^k\prod_{l\in P(T)}(x^{\mathfrak{d}(l)}-1)}{(x-1)^p}\times\underset{v\in S^T}{\prod}(x+\sum_{l\in L_v}\frac{x-1}{x^{\mathfrak{d}(l)}-1}).
\end{equation}
\textbf{Case I:} Let $T'$ be a uniform tree in $ \mathfrak{T}_{n,p}$. Then consider $T=T'$.\\
\\
\textbf{Case II:} Let $T'$ be a non-uniform tree. Let $T$ be a uniform tree in $\mathfrak{T}_{n,p}$ with seeding vertices $\{s_1,\ldots,s_k\}$, where each seeding vertex $s_j$ has $a_j$ seeds, for $1\leq j\leq k$. Let $l'$ and $l''$ be two seeds of a seeding vertex $v$ of $T$.

\noindent \textbf{Subcase I:} Suppose that both $l'$ and $l''$ are at distance $d$ from $v$. Assume that $l'$ is adjacent to a vertex $u$ in $T$. Delete the edge $ul'$ and add the edge $l'l''$. Let $T_1$ denote the resulting tree. Clearly, $u$ becomes a pendant vertex of $T_1$, and its distance from $v$ decreases to $d-1$. Moreover, $l''$ is no longer a pendant vertex of $T_1$, whereas $l'$ becomes a pendant vertex. Consequently, the distance of $l'$ from $v$ in $T_1$ is $d+1$. Thus, the connected zero forcing polynomial of $T_1$ is given by
\begin{equation}\label{equ.zerof.1}
    Z_C^{T_1}(x)=\frac{\underset{\substack{{l\in P(T)} \\ {l\neq l',
     l''}}}\prod(x^{\mathfrak{d}(l)}-1)}{x^{-k}(x-1)^p}\underset{\substack{v'\in S^T \\ v' \neq v}}{\prod}(x+\sum_{l\in L_v'}\frac{x-1}{x^{\mathfrak{d}(l)}-1})(x^{d-1}-1)(x^{d+1}-1)(x+\frac{x-1}{x^{d-1}-1}+\frac{x-1}{x^{d+1}-1}+\underset{\substack{l \in L_v, \\ l \neq l',l'' }}\sum \frac{x-1}{x^{\mathfrak{d}(l)}-1}).
\end{equation}
Let 
   $$A=\frac{x^k\underset{\substack{{l\in P(T)} \\ {l\neq l',
     l''}}}\prod(x^{\mathfrak{d}(l)}-1)}{(x-1)^p}\times\underset{\substack{v'\in S^T \\ v' \neq v}}{\prod}(x+\sum_{l\in L_{v'}}\frac{x-1}{x^{\mathfrak{d}(l)}-1}),$$ 
and $$B=\underset{\substack{l \in L_v, \\ l \neq l',l'' }}\sum \frac{x-1}{x^{\mathfrak{d}(l)}-1}.$$
So Equation (\ref{equ.zerof.1}) can be written as
\begin{equation}\label{equ.T}
    Z_C^{T_1}(x)= A \times (x^{d-1}-1)(x^{d+1}-1)(x+\frac{x-1}{x^{d-1}-1}+\frac{x-1}{x^{d+1}-1}+B).
\end{equation}
Similarly, Equation (\ref{equ.zerof}) can be written as
\begin{equation}\label{equ.T1}
Z_C^{T}(x)= A \times (x^{d}-1)^2(x+2\frac{x-1}{x^d-1}+B).  
\end{equation}
Using Equation \ref{equ.poly}
$$[(\sum_{i=0}^{d-1}x^i)^2x+2(\sum_{i=0}^{d-1}x^i)+B(\sum_{i=0}^{d-1}x^i)^2]$$
$$=[x(1+2x+3x^2+\cdots+(d-1)x^{d-2}+dx^{d-1}+(d-1)x^{d}+(d-2)x^{d+1}+\cdots+2x^{2d-3}+x^{2d-2})\times +2\sum_{i=0}^{d-1}x^i+B(\sum_{i=0}^{d-1}x^i)^2]$$
\begin{equation}\label{equ.compare1}
   =2+3x+4x^2+5x^3+6x^4\cdots+(d+1)x^{d-1}+dx^d+(d-1)x^{d+1}+\cdots+2x^{2d-2}+x^{2d-1}+B(\sum_{i=0}^{d-1}x^i)^2. 
\end{equation}

Again, using Equation \ref{equ.poly}, we have $$[(\sum_{i=0}^{d-2}x^i)(\sum_{i=0}^{d}x^i)x+\sum_{i=0}^{d-2}x^i+\sum_{i=0}^{d}x^i+B(\sum_{i=0}^{d-2}x^i)(\sum_{i=0}^{d}x^i)]$$
$$=[x(1+2x+\cdots+(d-1)x^{d-2}+(d-1)x^{d-1}+(d-1)x^{d}+\cdots+2x^{2d-3}+x^{2d-2}) +\sum_{i=0}^{d-2}x^i+\sum_{i=0}^{d}x^i+B(\sum_{i=0}^{d-2}x^i)(\sum_{i=0}^{d}x^i)]$$
\begin{equation}\label{equ.compare2}
 =2+3x+4x^2+5x^3+6x^4\cdots+(d+0)x^{d-1}+dx^d+(d-1)x^{d+1}+\cdots+2x^{2d-2}+x^{2d-1}+B(\sum_{i=0}^{d-2}x^i)(\sum_{i=0}^{d}x^i).   
\end{equation}

If we notice the Equations (\ref{equ.compare1}) and (\ref{equ.compare2}), then it is clear, $$(\sum_{i=0}^{d-2}x^i)(\sum_{i=0}^{d}x^i) \leq (\sum_{i=0}^{d-1}x^i)^2.$$
So,
\begin{equation}\label{equ.com}
   B(\sum_{i=0}^{d-2}x^i)(\sum_{i=0}^{d}x^i) \leq B(\sum_{i=0}^{d-1}x^i)^2. 
\end{equation}
Hence, using Equations (\ref{equ.compare1}), (\ref{equ.compare2}), and (\ref{equ.com}) from the Equations (\ref{equ.T}) and (\ref{equ.T1}), we can say 
$$Z_C^{T}(x) \geq Z_C^{T_1}(x).$$
\\
\textbf{Subcase II:} Suppose that the distances of the vertices $l'$ and $l''$ from $v$ are $d$ and $d+1$, respectively. Now, following the same construction as in Subcase I, we obtain a tree $T_1$ in which $u$ becomes a pendant vertex at distance $d-1$ from $v$, while $l'$ becomes a pendant vertex at distance $d+2$ from $v$. Thus, the connected zero forcing polynomial of $T_1$ is given by
\begin{equation}\label{equ.T2}
    Z_C^{T_1}(x)= A \times (x^{d-1}-1)(x^{d+2}-1)(x+\frac{x-1}{x^{d-1}-1}+\frac{x-1}{x^{d+2}-1}+B), 
\end{equation}
where the values of $A$ and $B$ are exactly the same as those in Subcase I. Again, following a similar argument and using Equation (\ref{equ.poly}), we obtain, from Equations (\ref{equ.T}) and (\ref{equ.T2}), that
$$Z_C^{T}(x) \geq Z_C^{T_1}(x).$$
\textbf{Subcase III:} Suppose that $v',v''\in \rm{st}(T)$, where $v'$ and $v''$ are the seeding vertices of $l'$ and $l''$, respectively. Assume that the distance of $l'$ from $v'$ is $d$, while the distance of $l''$ from $v''$ is $d+1$. Suppose that $l'$ is adjacent to a vertex $u$ in $T$. Delete the edge $ul'$ and add the edge $l'l''$. Let $T_1$ denote the resulting tree. Then
\begin{equation}\label{equ.3.}
    Z_C^{T_1}(x)= C \times (x^{d-1}-1)(x^{d+2}-1)(x+\frac{x-1}{x^{d-1}-1}+\underset{\substack{l \in L_{v'}, \\ l \neq l' }}\sum \frac{x-1}{(x^{\mathfrak{d}(l)}-1)})(x+\frac{x-1}{x^{d+1}-1}+\underset{\substack{l \in L_{v''}, \\ l \neq l'' }}\sum \frac{x-1}{x^{\mathfrak{d}(l)}-1}), 
\end{equation}
where 
\begin{equation}\label{equ.3}
   C=\frac{x^k\underset{\substack{{l\in P(T)} \\ {l\neq l',
     l''}}}\prod(x^{\mathfrak{d}(l)}-1)}{(x-1)^p}\times\underset{\substack{v\in S^T \\ v \neq v',v''}}{\prod}(x+\sum_{l\in L_{v}}\frac{x-1}{x^{\mathfrak{d}(l)}-1}). 
\end{equation}
Let
\begin{equation}\label{equ.3a}
    D_1=\underset{\substack{l \in L_{v'}, \\ l \neq l' }}\sum \frac{x-1}{(x^{\mathfrak{d}(l)}-1)},\ D_2=\underset{\substack{l \in L_{v'}, \\ l \neq l'' }}\sum \frac{x-1}{(x^{\mathfrak{d}(l)}-1)}.
\end{equation}
So, Equation (\ref{equ.3.}) can be written as 
\begin{equation}\label{equ.4}
       Z_C^{T_1}(x)= C \times (x^{d-1}-1)(x^{d+2}-1)(x+\frac{x-1}{x^{d-1}-1}+D_1)(x+\frac{x-1}{x^{d+2}-1}+D_2) 
\end{equation}
Similarly, Equation (\ref{equ.zerof}) can be written as
\begin{equation}\label{equ.5}
 Z_C^{T_1}(x)= C \times (x^{d}-1)(x^{d+1}-1)(x+\frac{x-1}{x^{d}-1}+D_1)(x+\frac{x-1}{x^{d+1}-1}+D_2)    
\end{equation}
Now, $$[x\sum_{i=0}^{d-2}x^i+D_1 \sum_{i=0}^{d-2}x^i] [x\sum_{i=0}^{d+1}x^i+D_2 \sum_{i=0}^{d+1}x^i]$$
\begin{equation}\label{equ.6}
  =x^2\sum_{i=0}^{d-2}x^i \sum_{i=0}^{d+1}x^i+x(D_1+D_2)\sum_{i=0}^{d-2}x^i \sum_{i=0}^{d+1}x^i+D_1D_2 \sum_{i=0}^{d-2}x^i \sum_{i=0}^{d+1}x^i.  
\end{equation}
Again, $$[x\sum_{i=0}^{d-1}x^i+D_1 \sum_{i=0}^{d-1}x^i] [x\sum_{i=0}^{d}x^i+D_2 \sum_{i=0}^{d}x^i]$$
\begin{equation}\label{equ.7}
=x^2\sum_{i=0}^{d-1}x^i \sum_{i=0}^{d}x^i+x(D_1+D_2)\sum_{i=0}^{d-1}x^i \sum_{i=0}^{d}x^i+D_1D_2 \sum_{i=0}^{d-1}x^i \sum_{i=0}^{d}x^i.    
\end{equation}
By above calculations and observing Equation (\ref{equ.poly}) it is clear that, if $r_1-m_1 \geq r_2-m_2$ then 
\begin{equation}\label{equ.lg}
  \sum_{i=0}^{m_1}x^i \sum_{j=0}^{r_1}x^j \leq \sum_{i=0}^{m_2}x^i \sum_{j=0}^{r_2}x^j.  
\end{equation}

So,
\begin{equation}\label{equ.8}
 \sum_{i=0}^{d-1}x^i \sum_{i=0}^{d}x^i \geq \sum_{i=0}^{d-2}x^i \sum_{i=0}^{d+1}x^i.   
\end{equation}
Hence, using Equations (\ref{equ.6}), (\ref{equ.7}), and (\ref{equ.8}), from Equations (\ref{equ.4}) and (\ref{equ.5}), we can clearly say
$$Z_C^{T}(x)\geq Z_C^{T_1}(x).$$
\textbf{Subcase IV:} Suppose that $v',v''\in \rm{st}(T)$, where $v'$ and $v''$ are the seeding vertices of $l'$ and $l''$, respectively. Assume that the distance of $l'$ from $v'$ and that of $l''$ from $v''$ are both equal to $d$. Suppose that $l'$ is adjacent to a vertex $u$ in $T$. Delete the edge $ul'$ and add the edge $l'l''$. Let $T_1$ denote the resulting tree. Then, following an argument similar to that in Subcase III, we obtain
 \begin{equation}
       Z_C^{T_1}(x)= C \times (x^{d-1}-1)(x^{d+1}-1)(x+\frac{x-1}{x^{d-1}-1}+D_1)(x+\frac{x-1}{x^{d+1}-1}+D_2), 
\end{equation}
\begin{equation}
 Z_C^{T_1}(x)= C \times (x^{d}-1)(x^{d}-1)(x+\frac{x-1}{x^{d}-1}+D_1)(x+\frac{x-1}{x^{d}-1}+D_2),    
\end{equation}
where the values of $C$, $D_1$, and $D_2$ are given in Equations (\ref{equ.3}) and (\ref{equ.3a}). By a similar way just like Subcase III, someone can prove
$$Z_C^{T}(x)\geq Z_C^{T_1}(x).$$
Now, if $T_1=T'$, then the proof is complete. Otherwise, we construct a tree $T_2$ from $T_1$ by following the same procedure used to construct $T_1$ from $T$. Clearly, we have $T_2\leq T_1$. Thus, by the preceding arguments and observing Equation (\ref{equ.lg}), we can inductively construct a sequence of trees terminating at $T'$ such that
$$
T'=T_s\leq T_{s-1}\leq \cdots \leq T_2\leq T_1\leq T.
$$
Hence, the proof is complete.
\end{proof}

\begin{corollary}
    Let $T\in \mathfrak{T}_{n,p}$, and $p_1$, $p_2$ be two pendant vertices of $T$. Let $p_1$ be adjacent to some vertex $v$ of $T$. Suppose distance of $p_1$ and $p_2$ from their respective seeding vertices (or vertex) be $d_1$ and $d_2$ respectively. Delete the edge between $v$ and $p_1$ and connect $p_1$ with $p_2$ by an edge. Let $T'$ be this new tree. Then $Z_C^{T}(x)\geq Z_C^{T'}(x)$ if and only if  $|d_1-d_2|\leq |(d_1-1)-(d_2+1)|$.
 \end{corollary}

\begin{corollary}\label{cor.starlike}
    Let $T$ and $ T' \in \mathfrak{T}_{n,p}$ be two non-isomorphic graphs with the vertices $\{v_1,\cdots,v_n\}$ which satisfy the following conditions.
\begin{enumerate}
\item $T$ and $T'$ are both starlike graphs.
\item  $T$ is a uniform tree.
\end{enumerate}
Then  $Z_C(T,i)\geq Z_C(T',i)$ for all $1\leq i\leq n$ and for any star graph $T'\in \mathfrak{T}_{n,p}$.
\end{corollary}

\begin{proof}
It is clear that, if we consider the set of all starlike graphs in $\mathfrak{T}_{n,p}$, then the set has a unique uniform tree. Hence, by Theorem \ref{thm.uniform}, above $T$ is the unique maximal element for the set of all starlike graphs in $\mathfrak{T}_{n,p}$.

\end{proof}

\begin{corollary}\label{thm 2.a}
 Let $T \in \mathfrak{T}_{n,3}$ be the uniform tree. Then $Z_C(T,i)\geq Z_C(T',i)$ for all $1\leq i\leq n$ and for any tree $T'\in \mathfrak{T}_{n,3}$.
\end{corollary}

\begin{lemma}\label{theorem.starlike}
Let $T\in \mathfrak{T}_{n,4}$. Then there exists a starlike graph $T'\in \mathfrak{T}_{n,4}$ such that $Z_C(T,i)\leq Z_C(T',i)$ for $1\leq i\leq n$.   
\end{lemma}
\begin{proof}
Let $T\in \mathfrak{T}_{n,4}$. If  $T$ is a starlike graph, then we can choose $T'=T$. Otherwise, $st(T)$ is a path with the two terminal vertices having external degree $2$, and all other vertices of $st(T)$ have external degree zero. Suppose, $T$ is not a graph of type $\mathcal{H}$. We now perform the following operation on $T$. Let $u_1, u_2,\ldots, u_l$ be the vertices of $st(T)$ with $u_i$ being adjacent to $u_{i+1}$. We do contraction on the edge between $u_1, u_2$ and add a new pendant vertex on a leaf having minimum distance from its seed $u_1$. We denote the obtained graph by $T_1$. Clearly $T_1\in \mathfrak{T}_{n,4}$. We will now show that $Z_C(T_1,i)\geq Z_C(T,i)$ for $1\leq i\leq n$. Let the distance of the two pendant vertices from their seeding vertex $u_1$ be $s_1,s_2$ respectively in $T$, with $s_1\leq s_2$. Then the connected zero forcing polynomial $Z_C^T(x)$ is of the form $$Z_C^T(x)=x^{l+2}(x(\frac{x^{s_1}-1}{x-1})(\frac{x^{s_2}-1}{x-1})+\frac{x^{s_1}-1}{x-1}+\frac{x^{s_2}-1}{x-1})g(x),$$ and the connected zero forcing polynomial $Z_C^{T_1}(x)$ is of the form $$Z_C^{T_1}(x)=x^{l+1}(x(\frac{x^{s_1+1}-1}{x-1})(\frac{x^{s_2}-1}{x-1})+\frac{x^{s_1+1}-1}{x-1}+\frac{x^{s_2}-1}{x-1})g(x),$$ where $g(x)$ is a polynomial with all coefficients non-negative. Now, $$x(\frac{x^{s_1+1}-1}{x-1})(\frac{x^{s_2}-1}{x-1})+\frac{x^{s_1+1}-1}{x-1}+\frac{x^{s_2}-1}{x-1}=x(1+x(\frac{x^{s_1}-1}{x-1}))(\frac{x^{s_2}-1}{x-1})+1+x(\frac{x^{s_1}-1}{x-1})+\frac{x^{s_2}-1}{x-1}.$$ So, $$Z_C^{T_1}(x)-Z_C^T(x)=x^{l+1}g(x)(1+\sum_{i=0}^{s_2-1}x^i)=x^{l+1}g(x)(1+\frac{x^{s_2}-1}{x-1}).$$ 
Thus, we have that $Z_C^{T_1}(x)-Z_C^T(x)$ is a polynomial with all coefficients non-negative. Hence, $Z_C(T_1,i)\geq Z_C(T,i)\  \mbox{for}\  1\leq i\leq n$. Continuing the above process, we can construct a tree $T_h \in \mathfrak{T}_{n,4}$ of type $\mathcal{H}$, such that $Z_C(\mathcal{H},i)\geq Z_C(T,i)\  \mbox{for}\ 1\leq i\leq n.$ We now show that for any tree $T_h$ of type $\mathcal{H}$, there exists a starlike graph $T_h'$ such that $$Z_C(T_h',i)\geq Z_C(T_h,i)\ \mbox{for}\ 1\leq i\leq n.$$ Let the distance of the four pendant vertices from their seeding vertices in the graph $T_h$ be $t_1,t_2$ and $t_3, t_4$, respectively, where the distances $t_1$ and $t_2$ correspond to the pendant vertices with the same seeding vertex, and the distances $t_3$ and $t_4$ correspond to the pendant vertices with the other seeding vertex. Also, assume that $t_1\leq t_2\leq t_3\leq t_4$. Construct the graph $T_h'$ by contracting the only edge in $st(T_h)$ and adding a new pendant vertex with the vertex having distance $t_1$. Then we have $$Z_C^{T_h}(x)=x^4(x(\frac{x^{t_1}-1}{x-1})(\frac{x^{t_2}-1}{x-1})+\frac{x^{t_1}-1}{x-1}+\frac{x^{t_2}-1}{x-1})(x(\frac{x^{t_3}-1}{x-1})(\frac{x^{t_4}-1}{x-1})+\frac{x^{t_3}-1}{x-1}+\frac{x^{t_4}-1}{x-1}),$$  $$Z_C^{T_h'}(x)=x^4(x(\frac{x^{t_1+1}-1}{x-1})(\frac{x^{t_2}-1}{x-1})(\frac{x^{t_3}-1}{x-1})(\frac{x^{t_4}-1}{x-1})+(\frac{x^{t_1+1}-1}{x-1})(\frac{x^{t_2}-1}{x-1})(\frac{x^{t_3}-1}{x-1})$$ $$+(\frac{x^{t_1+1}-1}{x-1})(\frac{x^{t_2}-1}{x-1})(\frac{x^{t_4}-1}{x-1})
+(\frac{x^{t_1+1}-1}{x-1})(\frac{x^{t_3}-1}{x-1})(\frac{x^{t_4}-1}{x-1})+(\frac{x^{t_2}-1}{x-1})(\frac{x^{t_3}-1}{x-1})(\frac{x^{t_4}-1}{x-1}))$$ $$=x^4((x(1+x(\frac{x^{t_1}-1}{x-1}))(\frac{x^{t_2}-1}{x-1})+(1+x(\frac{x^{t_1}-1}{x-1}))+(\frac{x^{t_2}-1}{x-1}))(\frac{x^{t_3}-1}{x-1})(\frac{x^{t_4}-1}{x-1})$$ $$+(1+x(\frac{x^{t_1}-1}{x-1}))(\frac{x^{t_2}-1}{x-1})(\frac{x^{t_4}-1}{x-1})+(1+x(\frac{x^{t_1}-1}{x-1}))(\frac{x^{t_2}-1}{x-1})(\frac{x^{t_3}-1}{x-1}))$$
$$=x^4(x(\frac{x^{t_3}-1}{x-1})(\frac{x^{t_4}-1}{x-1})(x(\frac{x^{t_1}-1}{x-1})(\frac{x^{t_2}-1}{x-1})+\frac{x^{t_1}-1}{x-1}+\frac{x^{t_2}-1}{x-1})+x(\frac{x^{t_1}-1}{x-1})(\frac{x^{t_2}-1}{x-1})$$ $$\times(\frac{x^{t_3}-1}{x-1}+\frac{x^{t_4}-1}{x-1})+(1+\frac{x^{t_2}-1}{x-1})(\frac{x^{t_3}-1}{x-1})(\frac{x^{t_4}-1}{x-1})+(\frac{x^{t_2}-1}{x-1})(\frac{x^{t_3}-1}{x-1}+\frac{x^{t_4}-1}{x-1})).$$
So, $$Z_C^{T_h'}(x)-Z_C^{T_h}(x)=x^4((1+\frac{x^{t_2}-1}{x-1})(\frac{x^{t_3}-1}{x-1})(\frac{x^{t_4}-1}{x-1})-(\frac{x^{t_1}-1}{x-1})(\frac{x^{t_3}-1}{x-1}+\frac{x^{t_4}-1}{x-1}))$$
$$=x^4((\frac{x^{t_3}-1}{x-1})(\frac{x^{t_4}-1}{x-1}-\frac{x^{t_1}-1}{x-1})+(\frac{x^{t_4}-1}{x-1})((\frac{x^{t_2}-1}{x-1})(\frac{x^{t_3}-1}{x-1})-\frac{x^{t_1}-1}{x-1}))$$
$$=x^4((\frac{x^{t_3}-1}{x-1})(\frac{x^{t_4}-x^{t_1}}{x-1})+(\frac{x^{t_4}-1}{x-1})(\frac{x^{t_3}-x^{t_1}}{x-1})+(\frac{x^{t_3}-1}{x-1})(\frac{x^{t_4}-1}{x-1})(\sum_{i=1}^{t_2-1}x^i)).$$

Since $t_1\leq t_3\leq t_4$, we find that the above polynomial has all coefficients non-negative. Hence, $Z_C(T_h',i)\geq Z_C(T_h,i)$ for $1\leq i\leq n$. This concludes the proof.
\end{proof}

\begin{theorem}\label{thm3.a}
 Let $T \in \mathfrak{T}_{n,4}$ be the uniform starlike tree. Then $Z_C(T,i)\geq Z_C(T',i)$ for all $1\leq i\leq n$ and for any graph $T'\in \mathfrak{T}_{n,4}$.
\end{theorem}

\begin{proof}
    By Corollary \ref{cor.starlike} and Lemma \ref{theorem.starlike}, the proof is trivial.
\end{proof}
Now, let us consider $P=5$. Then for $n=6$, $\mathcal{S}_{1,1,1,1,1,}$ is the unique tree in $\mathfrak{T}_{6,5}$. Again for $7\leq n\leq 8$, one can easily check that $T\in \mathcal{H}^n_{2,3}$ is the unique maximal element, where $T$ is the tree in Figure \ref{fig:result}.

\begin{center}
    
\begin{figure}
    \centering
    
\begin{tikzpicture}[scale=0.8,
    vertex/.style={circle, fill=black, inner sep=1.6pt},
    every node/.style={font=\large},
    line width=0.8pt
]

\begin{scope}[xshift=0cm,yshift=.5cm]

\draw
(0,0) -- (-0.8,0)
(0,0) -- (-0.7,0.65)
(0,0) -- (2.0,0)
(2.0,0) -- (2.75,0.7)
(2.0,0) -- (2.85,0)
(2.0,0) -- (2.85,-0.65);

\foreach \p in {
    (0,0),(-0.8,0),
    (-0.7,0.65),
    (2.0,0),
    (2.75,0.7),
    (2.85,0),
    (2.85,-0.65)
}
    \node[vertex] at \p {};
\node[below=0.9cm] at (1,0) {$T\in \mathcal{H}^7_{2,3}$};

\end{scope}

\begin{scope}[xshift=8cm,yshift=.5cm]

\draw
(0,0) -- (-0.8,0)
(0,0) -- (-0.7,0.65)
(0,0) -- (2.0,0)
(2.0,0) -- (2.75,0.7)
(2.0,0) -- (2.85,0) -- (3.45,0)
(2.0,0) -- (2.8,-0.65);

\foreach \p in {
    (0,0),(-0.8,0),
    (-0.7,0.65),
    (2.0,0),
    (2.75,0.7),
    (2.85,0),(3.45,0),
    (2.8,-0.65)
}
    \node[vertex] at \p {};

\node[below=0.9cm] at (1,0) {$T\in \mathcal{H}^8_{2,3}$};

\end{scope}
\end{tikzpicture}
\caption{Caption}
    \label{fig:result}

\end{figure}
\end{center}

 Now for $n\geq 9$, by the following \href{https://github.com/sanjaymath1729-ops/Connected-zero-forcing-polynomial-of-trees/blob/main/Finding%20maxima.py}{python program}, one can check that $\mathfrak{T}_{n,5}$ does not have a unique maximal element. In general someone can check the following results.
 \begin{result}\label{res.1}
     Let $T,\ T' \in \mathfrak{T}_{n,5}$, where $T'$ be the uniform starlike tree and $ T \in \mathcal{H}^n_{2,3}$ be a uniform tree. Also let $ \rm{st}(T)=\{v_1,\ v_2\}$ and degree of $v_1$, $v_2$ be $2$ and $3$ respectively. Let  $n=it+j(t+1)+2$, where $0\leq i,\ j\leq 5,\ i+j=5$, and $n\geq 9$. Then $T$ and $T'$ does not dominate by any other tree and they are never comparable when $T$ satisfies the following conditions.
 \begin{itemize}
        \item If $0\leq j \leq 3$, then distance of both the seeds of $v_1$ from $v_1$ is $t$.
        \item  If $4\leq j \leq 5$, then distance of every the seeds of $v_2$ from $v_2$ is $t+1$.
    \end{itemize}
  \begin{figure}[!ht]
\centering
\resizebox{0.5\textwidth}{!}{%
\begin{circuitikz}
\tikzstyle{every node}=[font=\normalsize]
\node at (5.5,7.5) [circ] {};
\draw (5.5,7.5) to[short] (4.75,7.5);
\draw (5.5,7.5) to[short] (5.5,6.75);
\draw (5.5,7.5) to[short] (5.5,8.25);
\draw [dashed] (5.5,8.25) -- (5.5,8);
\draw [dashed] (4.75,7.5) -- (4,7.5);
\draw [dashed] (5.5,6.75) -- (5.5,6.75);
\draw [dashed] (5.5,6.75) -- (5.5,6);
\draw [dashed] (1.25,7) -- (1.25,7);
\draw [short] (5.5,7.5) -- (6.5,7.5);
\draw [dashed] (6.5,7.5) -- (6.5,7.5);
\draw [dashed] (6.5,7.5) -- (7.25,7.5);
\node [font=\normalsize] at (5.5,9.25) {$h_1$};
\node [font=\normalsize] at (3.75,8.5) {$h_2$};
\node [font=\normalsize] at (2.75,7.5) {$h_3$};
\node [font=\normalsize] at (5.5,5.75) {$h_4$};
\node [font=\normalsize] at (8.5,7.5) {$h_5$};
\draw [short] (5.5,7.5) -- (5.5,7.5);
\draw [short] (5,8) -- (5,8);
\draw [short] (5.5,7.5) -- (4.5,8.25);
\draw [dashed] (3.75,8.75) -- (3,9.25);
\draw [dashed] (5.5,9.5) -- (5.5,10);
\draw [dashed] (1.5,7.5) -- (0.75,7.5);
\draw [dashed] (5.5,5.5) -- (5.5,4.75);
\draw [dashed] (9.75,7.5) -- (10.5,7.5);
\draw [short] (10.5,7.5) -- (11.5,7.5);
\draw [short] (5.5,10.25) -- (5.5,11);
\draw [short] (3,9.25) -- (2,10);
\draw [short] (0.75,7.5) -- (0,7.5);
\draw [short] (5.5,4.75) -- (5.5,4);
\node at (11.5,7.5) [circ] {};
\node at (5.5,11) [circ] {};
\node at (2,10) [circ] {};
\node at (0,7.5) [circ] {};
\node at (5.5,4) [circ] {};

\end{circuitikz}
}%
\caption{The uniform starlike tree $T'
\in \mathfrak{T}_{n,5}$.}
\label{fig.result1}
\end{figure}
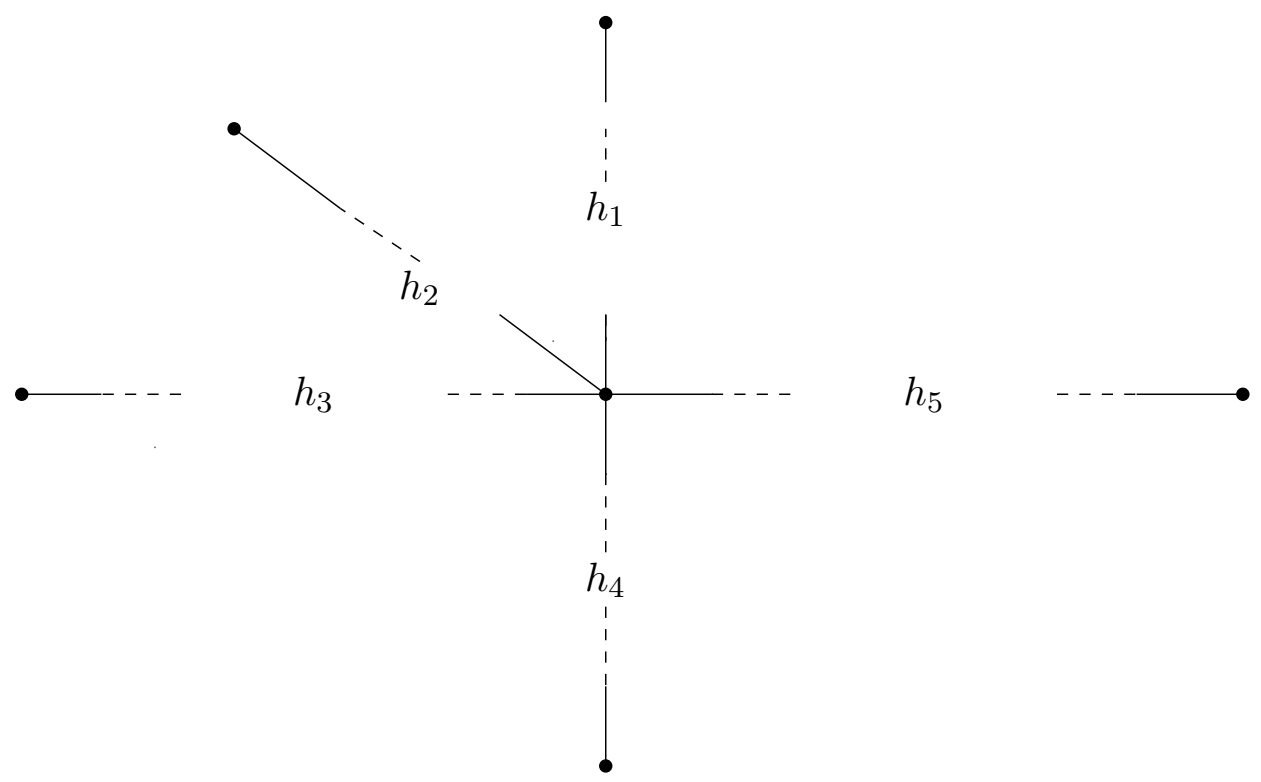
In Figure \ref{fig.result1}, for $1\leq i \leq 5$, $h_i$ is either $t$ or $t+1$ and $n=(h_1+\cdots+h_5)+1$.\par
In Figure \ref{fig:result.2}, for $0\leq j \leq 3$, $h_j$ is either $t$ or $t+1$, and $n=h_1+h_2+h_3+2t+2$. Again, for $4\leq j\leq 5$, $h_j$ is either $t$ or $t+1$, and $n=h_4+h_5+3(t+1)+2$.
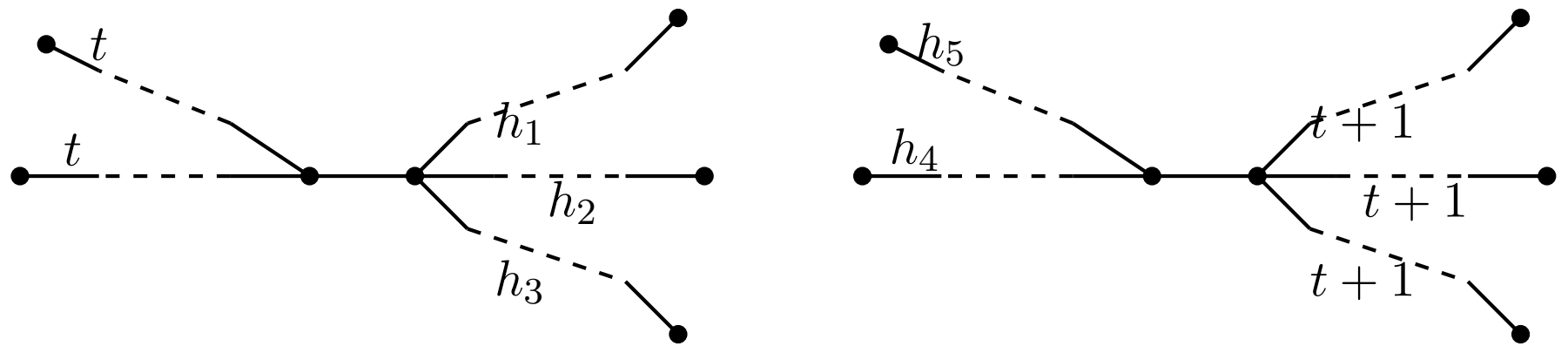
\begin{center}
    \begin{figure}
    \centering
    
\begin{tikzpicture}[scale=0.8,
    vertex/.style={circle, fill=black, inner sep=1.6pt},
    every node/.style={font=\large},
    line width=0.8pt
]

\begin{scope}[xshift=0cm,yshift=.5cm]

\node at (0,0) [circ] {};
\node at (1.0,0) [circ] {};
\draw (0,0) to[short] (1.0,0);
\draw (1,0) to[short] (1.5,.5);
\draw [dashed] (1.5,.5) -- (3,1);
\node [font=\normalsize] at (2,.5) {$h_1$};
\draw (3,1) to[short] (3.5,1.5);
\node at (3.5,1.5) [circ] {};

\draw (1,0) to[short] (1.75,0);
\draw [dashed] (1.75,0) -- (3,0);
\node [font=\normalsize] at (2.5,-.25) {$h_2$};
\draw (3,0) to[short] (3.75,0);
\node at (3.75,0) [circ] {};

\draw (1,0) to[short] (1.5,-.5);
\draw [dashed] (1.5,-.5) -- (3,-1);
\node [font=\normalsize] at (2,-1) {$h_3$};
\draw (3,-1) to[short] (3.5,-1.5);
\node at (3.5,-1.5) [circ] {};

\draw (0,0) to[short] (-.75,0);
\draw [dashed] (-0.75,0) -- (-2,0);
\node [font=\normalsize] at (-2.25,.25) {$t$};
\draw (-2,0) to[short] (-2.75,0);
\node at (-2.75,0) [circ] {};

\draw (0,0) to[short] (-.75,.5);
\draw [dashed] (-.75,.5) -- (-2,1);
\node [font=\normalsize] at (-2,1.25) {$t$};
\draw (-2,1) to[short] (-2.5,1.25);
\node at (-2.5,1.25) [circ] {};

\node[below=1.5cm] at (1,0) {$T\in \mathcal{H}^n_{2,3}$, when $0\leq j\leq 3$};

\end{scope}

\begin{scope}[xshift=8cm,yshift=.5cm]

\node at (0,0) [circ] {};
\node at (1.0,0) [circ] {};
\draw (0,0) to[short] (1.0,0);
\draw (1,0) to[short] (1.5,.5);
\draw [dashed] (1.5,.5) -- (3,1);
\node [font=\normalsize] at (2,.5) {$t+1$};
\draw (3,1) to[short] (3.5,1.5);
\node at (3.5,1.5) [circ] {};

\draw (1,0) to[short] (1.75,0);
\draw [dashed] (1.75,0) -- (3,0);
\node [font=\normalsize] at (2.5,-.25) {$t+1$};
\draw (3,0) to[short] (3.75,0);
\node at (3.75,0) [circ] {};

\draw (1,0) to[short] (1.5,-.5);
\draw [dashed] (1.5,-.5) -- (3,-1);
\node [font=\normalsize] at (2,-1) {$t+1$};
\draw (3,-1) to[short] (3.5,-1.5);
\node at (3.5,-1.5) [circ] {};

\draw (0,0) to[short] (-.75,0);
\draw [dashed] (-0.75,0) -- (-2,0);
\node [font=\normalsize] at (-2.25,.25) {$h_4$};
\draw (-2,0) to[short] (-2.75,0);
\node at (-2.75,0) [circ] {};

\draw (0,0) to[short] (-.75,.5);
\draw [dashed] (-.75,.5) -- (-2,1);
\node [font=\normalsize] at (-2,1.25) {$h_5$};
\draw (-2,1) to[short] (-2.5,1.25);
\node at (-2.5,1.25) [circ] {};

\node[below=1.5cm] at (1,0) {$T\in \mathcal{H}^n_{2,3}$, when $4\leq j\leq 5$};

\end{scope}
\end{tikzpicture}
\caption{}
    \label{fig:result.2}

\end{figure}
\end{center}

 \end{result}
\begin{result}\label{res.2}
 Let $T,\ T_1,\ T'\in \mathfrak{T}_{n,5}$, where $T$ and $T'$ are defined in  Result \ref{res.1}. Let st$(T_1)=\{v_1',v_2'\}$ and degree of $v_1'$ and $v_2'$ be $2$ and $3$ respectively. Also let $n=it+j(t+1)+2$, where $ i=2,\ j=3$, and $n\geq 9$. Then $T,\ T_1,\ T'\in \mathfrak{T}_{n,5}$ are not comparable by each other and none of them is dominated by any $T''\in \mathfrak{T}_{n,5}$ if $T'$ satisfies the following conditions.
 \begin{itemize}
    \item The distance of the seeds of $v_1$ from $v_1$ is $t$ and $t+1$ respectively.
        \item  The distance of the seeds of $v_2$ from $v_2$ is $t$, $t+1$, $t+1$ respectively.  
 \end{itemize}
 But it can be checked any tree in $\mathfrak{T}_{n,5}$ is dominated by atleast one of $T,\ T_1$ or $T'$.
\end{result}

Let us consider some examples in the support of the above results.
\begin{example}\label{example.p5}
Let $T,\ T_1,\ T'\in \mathfrak{T}_{n,5}$ be same as Results \ref{res.1} and \ref{res.2}.\\
\textbf{Case I:} Let $n=9$. Then \\
$Z_C^T(x)=6x^5+13x^6+11x^7+5x^8+x^9,\
Z_C^{T_1}(x)=6x^5+12x^6+11x^7+5x^8+x^9,$\\
$Z_C^{T'}(x)=5x^5+13x^6+12x^7+5x^8+x^9.$\\  
\textbf{Case II:} Let $n=10$. Then\\ 
  $Z_C^T(x)=6x^5+17x^6+19x^7+12x^8+5x^9+x^{10},\
Z_C^{T_1}(x)=6x^5+16x^6+19x^7+13x^8+5x^9+x^{10},$\\
$Z_C^{T'}(x)=5x^5+17x^6+22x^7+14x^8+5x^9+x^{10}.$\\ 
\textbf{Case III:} Let $n=11$. Then \\
  $Z_C^T(x)=6x^5+20x^6+29x^7+25x^8+14x^9+5x^{10}+x^{11},\
Z_C^{T_1}(x)=6x^5+19x^6+27x^7+24x^8+14x^9+5x^{10}+x^{11},$
$Z_C^{T'}(x)=5x^5+21x^6+35x^7+30x^8+15x^9+5x^{10}+x^{11}.$\\ 
\textbf{Case IV:} Let $n=12$. Then\\ 
$Z_C^T(x)=Z_C^{T_1}(x)=6x^5+23x^6+39x^7+41x^8+30x^9+15x^{10}+5x^{11}+x^{12},$\\
$Z_C^{T'}(x)=5x^5+21x^6+39x^7+43x^8+31x^9+15x^{10}+5x^{11}+x^{12}.$\\
\textbf{Case V:} Let $n=13$. Then\\ 
$Z_C^T(x)=6x^5+23x^6+43x^7+53x^8+47x^9+31x^{10}+15x^{11}+5x^{12}+x^{13},$\\
$Z_C^{T_1}(x)=6x^5+23x^6+42x^7+51x^8+46x^9+31x^{10}+15x^{11}+5x^{12}+x^{13},$\\
$Z_C^{T'}(x)=5x^5+21x^6+43x^7+56x^8+50x^9+32x^{10}+15x^{11}+5x^{12}+x^{13}.$\\
It is readily verified that, in every case, $T$ and $T'$ are incomparable. Furthermore, except in Case II, we have $T_1 \leq T$; in Case II, $T_1$ and $T$ are also incomparable. On the other hand, $T_1$ and $T'$ are incomparable in every cases. Consequently, in Case II, the three trees $T$, $T_1$, and $T'$ are pairwise incomparable.

Moreover, a direct computational verification using the \href{https://github.com/sanjaymath1729-ops/Connected-zero-forcing-polynomial-of-trees/blob/main/Finding%20maxima.py}{python program}, someone can check that all other trees in $\mathfrak{T}_{n,5}$ are dominated by at least one of $T$ or $T$ or $T'$, but there does not have any $T'' \in \mathfrak{T}_{n,5}$ which dominates both of $T$ and $T'$. 
\end{example}
\begin{observation}
  Since the patterns of $T$, $T_1$, and $T'$ are fixed, increasing the value of $n$ does not affect the comparability relations among $T$, $T_1$, and $T'$. Only the coefficients get increased in the exponents of $x$ in their connected zero forcing polynomials and whenever $n$ is increased by $5$, the connected zero forcing polynomial acquires five more additional power of $x$ relative to the corresponding polynomial for the preceding value of $n$.  
\end{observation}

\begin{example}\label{example.p8}
Let us consider $p=8$. Now, for $n=9$, the unique starlike tree $\mathcal{S}_{1,1,\ldots,1}$ is the only element in $\mathfrak{T}_{n,8}$. For $n=10,\ 11,\ 12$ using the \href{https://github.com/sanjaymath1729-ops/Connected-zero-forcing-polynomial-of-trees/blob/main/Finding%20maxima.py}{python program}, someone can check that the following $T_1$, $T_2$ and $T_3$ in Figure \ref{fig:p8.1} are the unique maximal tree in $\mathfrak{T}_{10,8}$, $\mathfrak{T}_{11,8}$ and $\mathfrak{T}_{12,8}$ respectively.\par
But for $n\geq 13$, there does not exist a unique maximal tree in $\mathfrak{T}_{n,8}$. For $n=13$, the following two maximal trees $T_1'$ and $T_2'$ in Figure \ref{fig:p8.1a} are not comparable to each,  as
$$Z_C^{T_1'}(x)=15x^8+44x^9+49x^{10}+27x^{11}+8x^{12}+x^{13},$$
$$Z_C^{T_2'}(x)=18x^8+45x^9+48x^{10}+27x^{11}+8x^{12}+x^{13}.$$
In fact for $n=17$, the following tress $T_4,\ T_5,\ T_6,\ T_7 \in \mathfrak{T}_{17,8}$ in Figure \ref{fig:placeholder} are maximal and they are not comparable to each other.
\begin{center}
    
\begin{figure}
    \centering
    
\begin{tikzpicture}[scale=0.8,
    vertex/.style={circle, fill=black, inner sep=1.6pt},
    every node/.style={font=\large},
    line width=0.8pt
]

\begin{scope}[xshift=0cm,yshift=4.5cm]

\draw
(0,0) -- (-0.8,0)
(0,0) -- (-0.7,0.65)
(0,0) -- (0,0.8)
(0,0) -- (-0.6,-0.7)
(0,0) -- (2.0,0)
(2.0,0) -- (2.75,0.7)
(2.0,0) -- (2.85,0)
(2.0,0) -- (2.85,-0.65)
(2.0,0) -- (2.0,-0.7);

\foreach \p in {
    (0,0),(-0.8,0),
    (-0.7,0.65),
    (0,0.8),
    (-0.6,-0.7),
    (2.0,0),
    (2.75,0.7),
    (2.85,0),
    (2.85,-0.65),
    (2.0,-0.7)
}
    \node[vertex] at \p {};
\node[below=0.9cm] at (1,0) {$T_1\in \mathfrak{T}_{10,8}$};

\end{scope}

\begin{scope}[xshift=8cm,yshift=4.5cm]

\draw
(0,0) -- (-0.8,0)
(0,0) -- (-0.7,0.65)
(0,0) -- (-0.6,-0.7)
(0,0) -- (2.0,0)
(2.0,0) -- (2.75,0.7)
(2.0,0) -- (2.85,0) -- (3.45,0)
(3.45,0) -- (4,0.7)
(3.45,0) -- (4,-0.7)
(3.45,0) -- (4.45,0)
(2.0,0) -- (2.8,-0.65);

\foreach \p in {
    (0,0),(-0.8,0),
    (-0.7,0.65),
    (-0.6,-0.7),
    (2.0,0),
    (2.75,0.7),
    (3.45,0),
     (4.45,0),
      (4,0.7),
       (4,-0.7),
    (2.8,-0.65)
}
    \node[vertex] at \p {};

\node[below=0.9cm] at (1,0) {$T_2\in \mathfrak{T}_{11,8}$};

\end{scope}

\begin{scope}[xshift=4cm,yshift=0.73cm]

\draw
(0,0) -- (-0.8,0)
(0,0) -- (-0.7,0.65)
(0,0) -- (-0.6,-0.7)
(0,0) -- (2.0,0)
(2.0,0) -- (2.75,0.7)
(2.0,0) -- (2.85,0) -- (3.45,0)
(3.45,0) -- (4,0.7)
(3.45,0) -- (4,-0.7)
(3.45,0) -- (4.45,0)
(-0.8,0) -- (-2,0)
(2.0,0) -- (2.8,-0.65);

\foreach \p in {
    (0,0),(-0.8,0),
    (-0.7,0.65),
    (-0.6,-0.7),
    (2.0,0),
    (2.75,0.7),
    (-2,0),
    (3.45,0),
     (4.45,0),
      (4,0.7),
       (4,-0.7),
    (2.8,-0.65)
}
    \node[vertex] at \p {};

\node[below=0.9cm] at (1,0) {$T_3\in \mathfrak{T}_{12,8}$};
\end{scope}
 \end{tikzpicture}
\caption{Unique maximal trees}
    \label{fig:p8.1}

\end{figure}
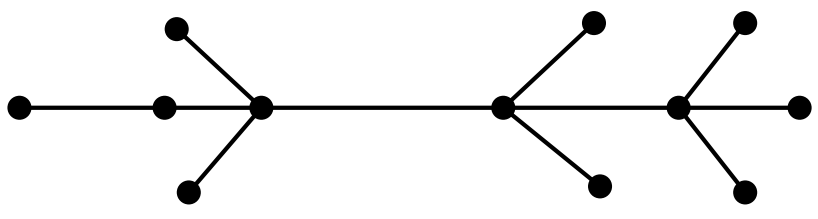
\end{center}

\begin{center}
    
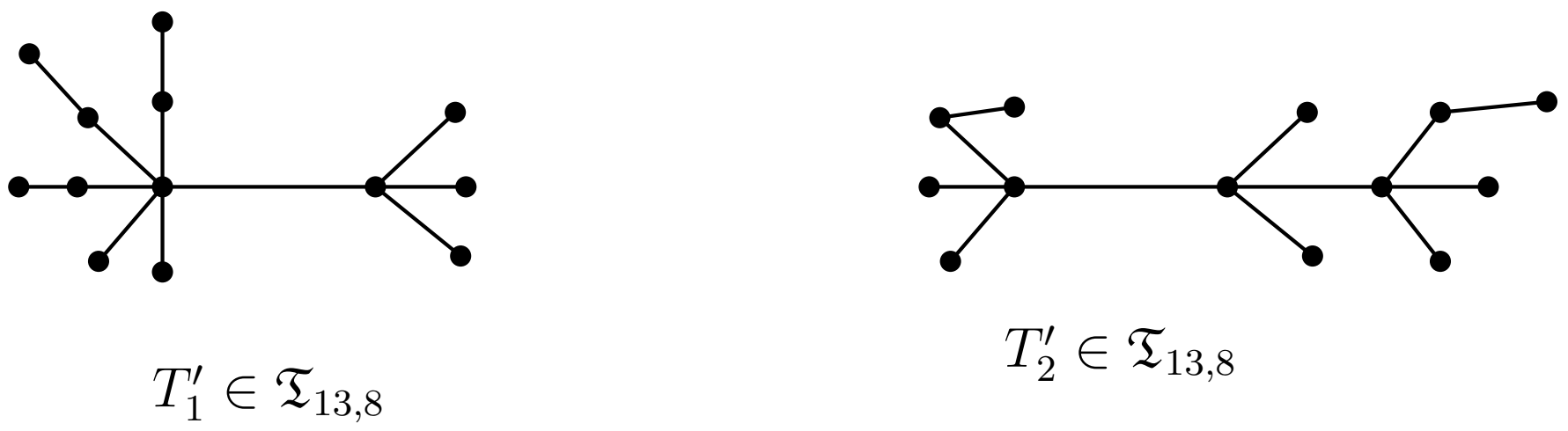
\begin{figure}
    \centering
    
\begin{tikzpicture}[scale=0.8,
    vertex/.style={circle, fill=black, inner sep=1.6pt},
    every node/.style={font=\large},
    line width=0.8pt
]

\begin{scope}[xshift=0cm,yshift=4.5cm]

\draw
(0,0) -- (-0.8,0) -- (-1.35,0)
(0,0) -- (-0.7,0.65) -- (-1.25,1.25)
(0,0) -- (-0.6,-0.7) 
(0,0) -- (0,0.8) -- (0,1.55)
(0,0) -- (0,-0.8) 
(0,0) -- (2.0,0)
(2.0,0) -- (2.75,0.7) 
(2.0,0) -- (2.85,0) 
(2.0,0) -- (2.8,-0.65);

\foreach \p in {
    (0,0),(-0.8,0),(-1.35,0),
    (-0.7,0.65),(-1.25,1.25),
    (-0.6,-0.7),
    (0,0.8),(0,1.55),
    (0,-0.8),
    (2.0,0),
    (2.75,0.7),
    (2.85,0),
    (2.8,-0.65)
}
    \node[vertex] at \p {};

\node[below=1.2cm] at (1,0) {$T_1'\in \mathfrak{T}_{13,8}$};
\end{scope}

\begin{scope}[xshift=8cm,yshift=4.5cm]

\draw
(0,0) -- (-0.8,0)
(0,0) -- (-0.7,0.65)
(-0.7,0.65) -- (0,.75)
(0,0) -- (-0.6,-0.7)
(0,0) -- (2.0,0)
(2.0,0) -- (2.75,0.7)
(2.0,0) -- (2.85,0) -- (3.45,0)
(3.45,0) -- (4,0.7)
(4,0.7) -- (5,0.8)
(3.45,0) -- (4,-0.7)
(3.45,0) -- (4.45,0)
(2.0,0) -- (2.8,-0.65);

\foreach \p in {
    (0,0),(-0.8,0),
    (0,.75),
    (-0.7,0.65),
    (-0.6,-0.7),
    (5,0.8),
    (2.0,0),
    (2.75,0.7),
    (3.45,0),
     (4.45,0),
      (4,0.7),
       (4,-0.7),
    (2.8,-0.65)
}
    \node[vertex] at \p {};

\node[below=0.9cm] at (1,0) {$T'_2\in \mathfrak{T}_{13,8}$};

\end{scope}

 \end{tikzpicture}
\caption{Non comparable maximal trees in $\mathfrak{T}_{13,8}$}
    \label{fig:p8.1a}

\end{figure}
\end{center}

\begin{center}
    
\begin{figure}
    \centering
    
\begin{tikzpicture}[scale=0.8,
    vertex/.style={circle, fill=black, inner sep=1.6pt},
    every node/.style={font=\large},
    line width=0.8pt
]

\begin{scope}[xshift=0cm,yshift=4.5cm]

\draw
(0,0) -- (-0.8,0) -- (-1.35,0)
(0,0) -- (-0.7,0.65) -- (-1.25,1.25)
(0,0) -- (0,0.8) -- (0,1.55)
(0,0) -- (-0.6,-0.7) -- (-1.2,-1.3)
(0,0) -- (2.0,0)
(2.0,0) -- (2.75,0.7) -- (3.35,1.3)
(2.0,0) -- (2.85,0) -- (3.45,0)
(2.0,0) -- (2.85,-0.65) -- (3.4,-1.15)
(2.0,0) -- (2.0,-0.7);

\foreach \p in {
    (0,0),(-0.8,0),(-1.35,0),
    (-0.7,0.65),(-1.25,1.25),
    (0,0.8),(0,1.55),
    (-0.6,-0.7),(-1.2,-1.3),
    (2.0,0),
    (2.75,0.7),(3.35,1.3),
    (2.85,0),(3.45,0),
    (2.85,-0.65),(3.4,-1.15),
    (2.0,-0.7)
}
    \node[vertex] at \p {};

\node[below=0.9cm] at (1,0) {$T_4$};

\end{scope}

\begin{scope}[xshift=8cm,yshift=4.5cm]

\draw
(0,0) -- (-0.8,0) -- (-1.35,0)
(0,0) -- (-0.7,0.65) -- (-1.25,1.25)
(0,0) -- (-0.6,-0.7) -- (-1.2,-1.3)
(0,0) -- (0,0.8) -- (0,1.55)
(0,0) -- (0,-0.8) -- (0,-1.55)
(0,0) -- (2.0,0)
(2.0,0) -- (2.75,0.7) -- (3.35,1.3)
(2.0,0) -- (2.85,0) -- (3.45,0)
(2.0,0) -- (2.8,-0.65);

\foreach \p in {
    (0,0),(-0.8,0),(-1.35,0),
    (-0.7,0.65),(-1.25,1.25),
    (-0.6,-0.7),(-1.2,-1.3),
    (0,0.8),(0,1.55),
    (0,-0.8),(0,-1.55),
    (2.0,0),
    (2.75,0.7),(3.35,1.3),
    (2.85,0),(3.45,0),
    (2.8,-0.65)
}
    \node[vertex] at \p {};

\node[below=0.9cm] at (1,0) {$T_5$};

\end{scope}

\begin{scope}[xshift=0cm,yshift=0.73cm]

\draw
(0,0) -- (-0.8,0) -- (-1.35,0)
(0,0) -- (-0.7,0.65) -- (-1.25,1.25)
(0,0) -- (-0.6,-0.7) -- (-1.2,-1.3)
(0,0) -- (2.0,0)
(2.0,0) -- (1.55,-0.8)
(2.0,0) -- (2.45,-0.8)
(2.0,0) -- (4.0,0)
(4.0,0) -- (4.7,0.7) -- (5.3,1.3)
(4.0,0) -- (4.85,0) -- (5.45,0)
(4.0,0) -- (4.7,-0.65) -- (5.3,-1.15);

\foreach \p in {
    (0,0),(-0.8,0),(-1.35,0),
    (-0.7,0.65),(-1.25,1.25),
    (-0.6,-0.7),(-1.2,-1.3),
    (2.0,0),(1.55,-0.8),(2.45,-0.8),
    (4.0,0),
    (4.7,0.7),(5.3,1.3),
    (4.85,0),(5.45,0),
    (4.7,-0.65),(5.3,-1.15)
}
    \node[vertex] at \p {};

\node[below=1.53cm] at (1.1,0) {$T_6$};

\end{scope}

\begin{scope}[xshift=9.15cm,yshift=0.73cm]

\draw
(0,0) -- (-0.8,0) -- (-1.4,0)
(0,0) -- (-0.7,0.7) -- (-1.3,1.3)
(0,0) -- (0,0.8) -- (0,1.5)
(0,0) -- (0.7,0.7) -- (1.3,1.3)
(0,0) -- (0.8,0) -- (1.4,0)
(0,0) -- (0.7,-0.7) -- (1.3,-1.3)
(0,0) -- (0,-0.8) -- (0,-1.5)
(0,0) -- (-0.7,-0.7) -- (-1.3,-1.3);

\foreach \p in {
    (0,0),
    (-0.8,0),(-1.4,0),
    (-0.7,0.7),(-1.3,1.3),
    (0,0.8),(0,1.5),
    (0.7,0.7),(1.3,1.3),
    (0.8,0),(1.4,0),
    (0.7,-0.7),(1.3,-1.3),
    (0,-0.8),(0,-1.5),
    (-0.7,-0.7),(-1.3,-1.3)
}
    \node[vertex] at \p {};

\node[below=1.6cm] at (0,0) {$T_7$};

\end{scope}

\end{tikzpicture}
\caption{Non comparable maximal trees in $\mathfrak{T}_{17,8}$}
    \label{fig:placeholder}
\end{figure}
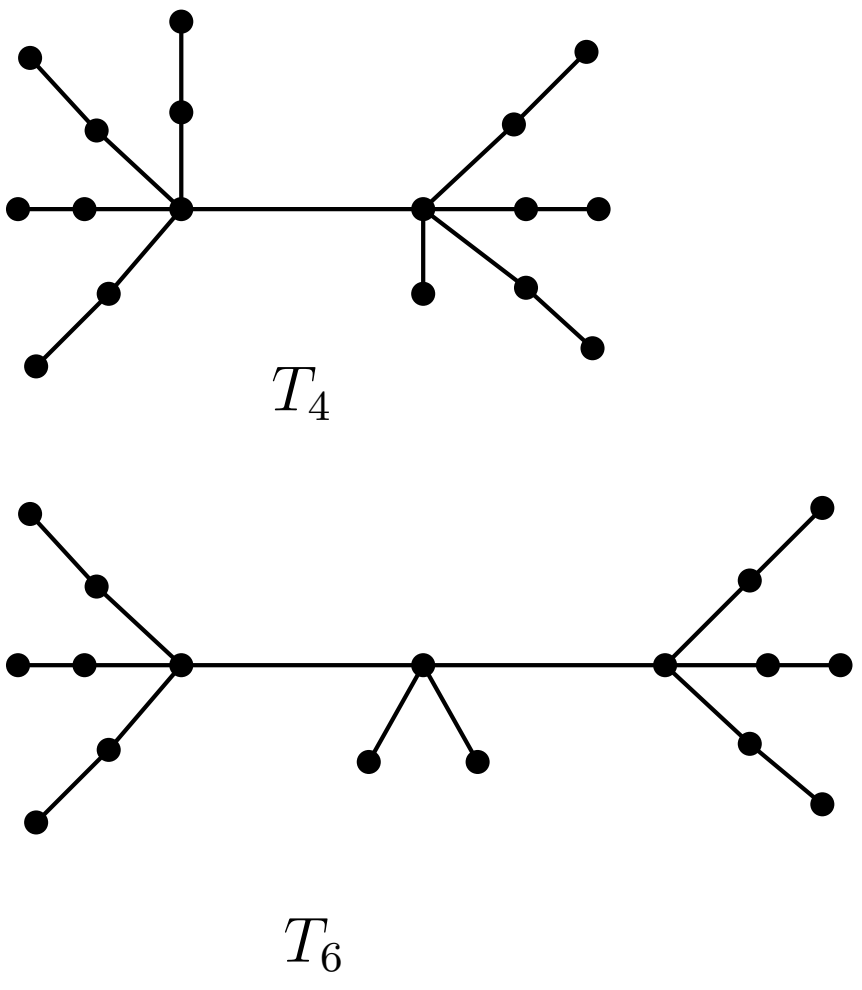
\end{center}

 By the \href{https://github.com/sanjaymath1729-ops/Connected-zero-forcing-polynomial-of-trees/blob/main/Computing%20the%20connected%20zero%20forcing%20polynomial%20of%20trees.py}{python program}, someone can check
 \begin{align*}
   Z_C^{T_4}(x)=16x^8+92x^9+230x^{10}+333x^{11}+316x^{12}+211x^{13}+102x^{14}+35x^{15}+8x^{16}+x^{17},\\ 
   Z_C^{T_5}(x)=15x^8+88x^9+225x^{10}+333x^{11}+321x^{12}+215x^{13}+103x^{14}+35x^{15}+8x^{16}+x^{17},\\
   Z_C^{T_6}(x)=18x^8+93x^9+212x^{10}+289x^{11}+270x^{12}+184x^{13}+92x^{14}+33x^{15}+8x^{16}+x^{17},\\
   Z_C^{T_7}(x)=8x^8+57x^9+176x^{10}+308x^{11}+336x^{12}+238x^{13}+112x^{14}+36x^{15}+8x^{16}+x^{17}.
 \end{align*}
It is clear that none of $T_4$, $T_5$, $T_6$ and $T_7$ are comparable to each other. Also, by the \href{https://github.com/sanjaymath1729-ops/Connected-zero-forcing-polynomial-of-trees/blob/main/Finding%20maxima.py}{python program}, someone can check none of them is dominated by any $T\in \mathfrak{T}_{17,8}$.
\end{example}
\begin{observation}
    From Example \ref{example.p5}, we know, for $p=5$, till $n=8$ we have a unique maximal tree in $\mathfrak{T}_{n,5}$. But for $n\geq 8$, we do not have a unique maximal tree. In fact we found $3$ maximal tree which are not comparable to each other and all other trees in $\mathfrak{T}_{n,5}$ are dominated by atleast one of them.\par
    Again, from Example \ref{example.p8} we know, for $p=8$, till $n=12$, there exists a unique maximal tree. But for $n\geq 13$, we do not have a unique maximal tree. In fact we found $4$ maximal trees in $\mathfrak{T}_{17,8}$ which are not comparable to each other.
\end{observation}


\section{Concluding Remarks and Open Problems}

Zero forcing sets and connected zero forcing sets have applications in various areas of mathematics and related disciplines. Although zero forcing sets have been extensively studied, the theory of connected zero forcing sets remains relatively unexplored. In this paper, we investigate several properties of connected zero forcing sets of trees by introducing the connected zero forcing polynomial $Z_C^T(x)$ of a tree $T$. We prove that this polynomial is always log-concave. Furthermore, for fixed values of $n$ and $p$, we determine the minimal elements in $\mathfrak{T}_{n,p}$. Through several illustrative examples, we also demonstrate that there need not exist a unique maximal element in $\mathfrak{T}_{n,p}$. This observation motivates the following conjecture.

\begin{conjecture}
For $p\geq 5$, there exists a sufficiently large integer $N$ (depending on $p$) such that, for every $n\geq N$, there is no unique tree $T\in\mathfrak{T}_{n,p}$ satisfying
$$
T\geq T',\qquad \text{for every }T'\in\mathfrak{T}_{n,p}.
$$
More generally, as the number of pendant vertices increases, the number of pairwise incomparable maximal elements in $\mathfrak{T}_{n,p}$ tends to increase for sufficiently large values of $n$, where the threshold may depend on $p$.
\end{conjecture}

In this paper, we investigate several properties of connected zero forcing polynomials of finite trees. Our results reveal that these polynomials exhibit properties that can differ significantly from those of the usual zero forcing polynomials of trees. Since unicyclic graphs are among the graph classes most closely related to trees, it would be interesting to investigate the connected zero forcing polynomials of unicyclic graphs. In particular, we propose the following problem.

\begin{problem}
Determine the connected zero forcing polynomial of a unicyclic graph. Is this polynomial always log-concave? Furthermore, for fixed positive integers $n$ and $p$, do there exist unicyclic graphs $U_{\min}$ and $U_{\max}$, each having $n$ vertices and $p$ pendant vertices, such that
$$
U_{\min}\leq U\leq U_{\max}
$$
for every unicyclic graph $U$ with $n$ vertices and $p$ pendant vertices?
\end{problem}

Throughout our investigation of the minimal and maximal elements of trees, we have fixed the number of pendant vertices. It is natural to ask what happens when this restriction is removed. More precisely, consider two trees $T_1\in\mathfrak{T}_{n,p_1}$ and $T_2\in\mathfrak{T}_{n,p_2}$, where $p_1$ and $p_2$ need not be equal. Let $Z_C^{T_1}(x)$ and $Z_C^{T_2}(x)$ denote their respective connected zero forcing polynomials. This leads to the following problem.

\begin{problem}
Given two trees $T_1\in\mathfrak{T}_{n,p_1}$ and $T_2\in\mathfrak{T}_{n,p_2}$, can we always find an integer $l$, with $1\leq l\leq n$, such that
$$
Z_C(T_1,i)\leq Z_C(T_2,i)
\quad\text{for }1\leq i\leq l
$$
and
$$
Z_C(T_1,i)\geq Z_C(T_2,i)
\quad\text{for }l+1\leq i\leq n?
$$
\end{problem}

Finally, although this paper focuses on trees, it is natural to extend the investigation to arbitrary finite graphs. Such an extension could provide a broader understanding of the structural and algebraic properties of connected zero forcing sets. In this direction, we propose the following problem.

\begin{problem}
Let $G$ be an arbitrary finite graph. Determine its connected zero forcing polynomial $Z_C^G(x)$. Is this polynomial always log-concave? Furthermore, given two finite graphs $G_1$ and $G_2$, each having $n$ vertices and $p$ pendant vertices, what relationships, if any, exist between their connected zero forcing sets and their corresponding connected zero forcing polynomials? More generally, what can be said about these relationships when the two graphs have different numbers of pendant vertices?
\end{problem}
There is a close relation between Restricted power domination and zero forcing number \cite{bozem}. In \cite{bur}, Burgarth et at investigated how the graph theoretic concept of a zero forcing set impacts the controllability property. In \cite{curt} and \cite{english}, the authors worked on zero forcing on random graphs. In \cite{group} and \cite{huang}, the authors used the zero forcing number to bound the minimum rank or  maximum nullity of a graph. In \cite{goldberg}, Goldberg and Berman provided a bound for the maximum nullity of a matrix with a given sign pattern.



\section*{Acknowledgments}

We thank Dr.~Anurag Singh for several conversations on this work. We also thank Dr.~Boris Brimkov for bringing  article \cite{Newton} and the open question posed there, in our attention. The first-named author thanks Harish-Chandra Research Institute, India, for providing funding (in the form of post-doctoral fellowship) during this work. The second-named author expresses gratitude to the Department of Mathematical Sciences IISER Mohali, India for providing funding (in the form of post-doctoral fellowship) during this work.

\section*{AI Disclosure}
The authors acknowledge the assistance provided by the LLM Claude (Anthropic; model Claude Opus~5.5), while figuring out some results in section \ref{Sec: Log}. All other results provided in any other section are found as well as written from scratch entirely by the authors. In addition, all mathematical statements and proofs are independently written and verified by the authors.

\section*{Conflict of interest}

The authors declare that there are no conflict of interests.

\bibliographystyle{plain}

\end{document}